\documentclass[10pt]{amsart}

\usepackage{times}
\usepackage{amsmath,  amssymb,slashed,url,bm,upgreek,amsthm,esint}
\usepackage{graphicx,enumerate}
\usepackage{color}
\usepackage{fancybox}
\usepackage{hyperref}
\newcounter{intro}

\newtheorem{thm}{Theorem}[section]
\newtheorem{lem}[thm]{Lemma}
\newtheorem{prop}[thm]{Proposition}
\newtheorem{cor}[thm]{Corollary}
\newtheorem{defi}[thm]{Definition}
\newtheorem{rem}[thm]{Remark}

\newcommand{\cref}[1]{Corollary~\ref{#1}}

\newcommand{\lref}[1]{Lemma~\ref{#1}}
\newcommand{\pref}[1]{Proposition~\ref{#1}}

\DeclareMathOperator{\un}{\mathbf{1}}
\DeclareMathOperator{\eucl}{eucl}

\DeclareMathOperator{\supp}{supp}

\DeclareMathOperator{\scal}{Scal}
\DeclareMathOperator{\diam}{diam}

\DeclareMathOperator{\vol}{Vol}
\def\R{\mathbb R}\def\N{\mathbb N}\def\Z{\mathbb Z}\def\bS{\mathbb S}\def\bB{\mathbb B}\def\bT{\mathbb T}
\def\cB{\mathcal B}
\def\cC{\mathcal C}

\def\cM{\mathcal M}
\def\loc{\mathrm{loc}}\DeclareMathOperator{\Q}{Q}

\def\ra{\rangle}
\def\la{\langle}

\newcommand{\dv}{d{\rm \mu}} 
\newcommand{\norm}[1]{\left|\!\left|#1\right|\!\right|}
\newcommand{\abs}[1]{\left|#1\right|}
\newcommand{\ds}{\displaystyle}
\newcommand{\e}{\varepsilon}
\renewcommand{\epsilon}{\varepsilon}

\author{Clara L. Aldana\\ University of Luxembourg\and 
 Gilles Carron  \\ Universit\'e de Nantes  \and 
Samuel Tapie \\ Universit\'e de Nantes}

\email{clara.aldana\@@posteo.net}
\email{gilles.carron\@@univ-nantes.fr}
\email{samuel.tapie\@@univ-nantes.fr}
\title{$A_\infty$ weights and compactness of conformal metrics under $L^{n/2}$ curvature bounds}

\begin{document}

\begin{abstract}

We study sequences of conformal deformations of a smooth closed Riemannian manifold of dimension $n$, assuming uniform volume bounds and $L^{n/2}$ bounds on their scalar curvatures. Singularities may appear in the limit. Nevertheless, we show that under such bounds the underlying metric spaces are pre-compact in the Gromov-Hausdorff topology. Our study is based on the use of \emph{$A_\infty$-weights} from harmonic analysis, and provides geometric controls on the limit spaces thus obtained. Our techniques also show that any conformal deformation of the Euclidean metric on $\mathbb R^n$ with infinite volume and finite $L^{n/2}$ norm of the scalar curvature satisfies the Euclidean isoperimetric inequality.

\end{abstract}

\maketitle

\markboth{$A_\infty$ weights and compactness of conformal metrics}{C. Aldana, G. Carron, S. Tapie}

\section{Introduction}

This paper is devoted to the question of sequential compactness of Riemannian metrics inside a conformal class given uniform bounds on their volume and integral bounds on their curvatures. It was shown by Gursky in \cite{Gursky} that (unless $g_0$ is conformal to the round sphere), if the Riemann curvature tensor is bounded in $L^p$ for some $p>\frac n 2$ and if the volume is uniformly controlled, then the sequence of metric is pre-compact in the $\cC^\alpha$ topology for some $\alpha\in (0,1)$. 

It was later shown that no such compactness can hold when the curvature tensor is bounded in $L^{n/2}$. Indeed, in \cite{CGW} Chang, Gursky and Wolff constructed two examples of families of sequences of conformal metrics with uniform volume and $L^{n/2}$ curvature bounds which have no convergent sub-sequences. Let us briefly recall these examples.

\subsection*{Example 1: spherical blowup}
Let $(\bS^n, g_S)$ be the standard sphere, $N\in \bS^n$ be some fixed point and $S$ be its antipodal pointÂ§. Let $\sigma: \bS^n\backslash \{N\} \to\R^n$ be the stereographic projection with pole $N$, and let $\mathrm{eucl}$ denote the Euclidean metric on $\R^n$. For all $\lambda>0$, let $g_\lambda = \sigma^*(\lambda \ \mathrm{eucl})$: it is well known that $g_\lambda$ extends to a smooth metric on $\bS^n$, conformally equivalent to $g_S$. 

Note that since $g_\lambda$ is merely a pullback of $g_S$ by some conformal diffeomorphisms, the volume of $(\bS^n, g_\lambda)$ and all $L^p$ norms of the curvature tensor are preserved. Moreover, as $\lambda\to +\infty$, the Riemannian measures ${\dv}_{g_\lambda}$ converge to a Dirac mass at $S$. 

It should also be noticed that for any $R>0$, the $g_S$-ball $\bB(S, R)$ of radius $R$ contains asymptotically all the curvature of $g_\lambda$. In particular, we have
$$
\lim_{\lambda\to +\infty} \left(\int_{\bB(S, R)} \left(\scal_{g_\lambda}\right)^{\frac n2} \dv_{g_\lambda}\right)^{\frac 2n} = \upalpha(n,2),
$$
where
$\upalpha(n,2)=\sigma_n^{\frac 2n} n(n-1)$ is the $L^{n/2}$-total scalar curvature of the standard round unit sphere, and $\sigma_n = \frac{2\pi^{n/2}}{\Gamma(n/2)}$ is its $n$-th volume. 

\medskip
It was shown in \cite{CGW} that this family of conformal deformations of the sphere can be glued on any compact Riemannian manifold $(M, g_0)$ at any point $x_0$, giving rise to a sequence of conformal deformations $(g_\lambda = e^{2 f_\lambda} g_0)_{\lambda>0}$ such that:
\begin{itemize}
\item the volume and diameter of the complement of any open ball centered at $x_0$ go to $0$ as $\lambda\to +\infty$;
\item the $L^{n/2}$ norm of the Riemann curvature tensor are uniformly bounded;
\item the sequence of Riemannian measures $d{\rm \mu}_{g_\lambda}$ develop a Dirac mass at $x_0$ and the volume stays in some fixed interval $[v,V]$. 
\end{itemize}
In particular, this implies that the associated sequence of distances $(d_{g_\lambda})_{\lambda>0}$ is not precompact for the uniform $\cC^0$ topology. Once again, in this example, for any fixed $R>0$, the $g_0$-ball $\bB(x_0, R)$ of radius $R$ contains asymptotically a large part of the scalar curvature of $g_\lambda$. In particular it can be shown that, in this case, we have
$$
\liminf_{\lambda \to +\infty} \left(\int_{\bB(x_0, R)} \left(\scal^+_{g_\lambda}\right)^{\frac n2}\dv_{g_\lambda}\right)^{\frac 2n} \geq \upalpha(n,2).
$$

We will see below that if we consider sequences of conformal measures whose volume and $L^{n/2}$-total \emph{scalar} curvature are uniformly bounded, and the scalar curvature concentration (in $L^{n/2}$ norm) in a $g_0$-ball is \emph{strictly less} than $\upalpha(n,2)$, then this kind of spherical blowup cannot occur.

\subsection*{Example 2: Schr\"odinger type blowup}
The second construction in \cite{CGW} exploits the Schr\"odinger-type structure of the conformal Laplacian, where the potential is only in $L^{n/2}$. We present a simple example inspired by their construction.

Let $(M, g_0)$ be a compact smooth Riemannian manifold, let $x_0\in M$ and $R_0$ be such that $\text{inj}_{x_0}(M, g_0) \geq 2 R_0$, where  $\text{inj}_{x_0}(M, g_0)$ denotes the injectivity radius of $g_0$ at $x_0$. For $x\in M$, define
\begin{equation}\label{eq:SchrodEx}
f_\infty(x) = \left\{\begin{array}{ccc}
					\sqrt{\left|\ln \left(\frac{d(x_0, x)}{R_0}\right)\right|} & {\rm if} & 0<d(x_0, x)\leq R_0\\
					0 & {\rm if} & d(x_0, x)\geq 2R_0
				\end{array}\right.
\end{equation}
with uniformly bounded second derivatives when $R_0\leq d(x_0, x)\leq 2R_0$. Let $(f_k : M\to \R)_{k\in \N}$ be a sequence of smooth maps which converges to $f_\infty$ in $W^{2, \frac n 2}$. We consider for all $k\geq 0$ the conformal metric $g_k = e^{2f_k} g_0$. The volume of $(M, g_k)$ is uniformly bounded, and it follows from Theorem 1.159 in \cite[p 58]{Besse} that the $L^{\frac n 2}$ norm of the Riemannian curvature of $(M, g_{k})$ are uniformly bounded. Note that eventhough the sequence of Riemannian metrics $(g_k)$ will blow up at $x_0$, the diameter of $(M, g_k)$ remains bounded and the associated volume measures will not develop a Dirac mass. Using the sequence of metrics $(\tilde g_k = e^{-2f_k} g_0)_{k\in \N}$, we get an analogous example where the metric vanishes at $x_0$ as $k\to +\infty$. 

Analogously, one can construct examples of sequences of conformal deformations whose volumes and $L^{\frac n 2}$ norms of the scalar curvatures are uniformly bounded, for which the limit metric blows up (or goes to $0$) at any countable (even dense) set of points.

\medskip
The reader may easily show that in this example, eventhough the sequence of conformal factors is unbounded, the sequence of associated distances $(d_{g_k})_{k\in \N}$ converges for the $\mathcal C^0$-uniform topology to the Riemannian distance associated to the (singular) metric $g_\infty = e^{2 f_\infty} g_0$ on $M$. 

\medskip

Our first theorem is a precompactness result which follows from an elementary application of Sobolev embeddings.

\begin{thm}\label{GHLn/2}
Let $(M,g_0)$ be a compact smooth Riemannian manifold and let $\cM_{V}^{\delta}$ be the space of conformal smooth metrics on $M$, $g_f = e^{2f} g_0$ with associated Riemannian measure $d\mu_f=e^{nf}\mathrm{dvol}_{g_0}$, such that there exists $V>0$, $\delta>0$ such that: 
\begin{enumerate}
\item $\ds \mu_f(M) \leq V$, and
\item $\ds \norm{(\scal_{g_f})_+}_{L^{\frac{n}{2}}(M, \dv_{f})}\le\upalpha(n,2)-\delta$.
\end{enumerate}
Then for some $\alpha\in (0,1)$, the set of distances $\{d_{g_f} \; ; \; g_f\in \cM_{V}^{\delta}\}$ is precompact for the $\cC^\alpha$ topology on $M\times M$. In particular, the set of metric spaces $\{(M,d_{g_f}), g_f\in \cM_{V}^{\delta}\}$ is precompact for the Gromov-Hausdorff topology.
\end{thm}
It follows from the Example 1 presented above that the upper bound 
$$
\norm{(\scal_{g_f})_+}_{L^{\frac{n}{2}}(M, \dv_{f})}\le\upalpha(n,2)-\delta = \sigma_n^{\frac 2n} n(n-1)-\delta
$$ is optimal to get $\cC^\alpha$ precompactness of the sequence of distances.

Theorem \ref{GHLn/2} does not exclude the possible collapsing of some part of the manifold in the  limit of such conformal deformations. We will show in Corollary \ref{Gurskyimprove} that adding to the previous hypotheses a lower bound on the volume and a $L^p$ bound on the \emph{scalar} curvature for some $p>\frac n 2$ is enough to recover $\cC^{\alpha}$ precompactness of the sequence of Riemannian metrics as in \cite{Gursky}.

The rest of our paper is devoted to the exploration of conformal deformations for which only volume bounds and $L^{\frac n 2}$ bounds on the scalar curvature are imposed. We will see in particular that, no collapsing can occur if the $L^{n/2}$ norm of the scalar curvature is small at small scales.

\medskip

Our ideas come from a paper by Yi Wang \cite{Wang}, who was interested in conformal deformations with integral bounds on the $Q$-curvature related to the so-called Muckenhoupt weights, or $A_\infty$-weights. If $(M,g)$ is a closed manifold, a non-negative $L^1_{loc}$ function $w$ is called an \emph{$A_\infty$-weight with respect to $g$} if there exist $C>0, q>1$ such that for all $g$-geodesic balls $\bB\subset M$,
$$\left(\fint_{\bB} w^qd\mu_g\right)^{\frac{1}{ q}}\le C \fint_{\bB} wd\mu_g,$$
where $\mu_g$ is the volume measure associated to $g$.
Several equivalent characterizations of these $A_\infty$-weights will be presented in Section \ref{sec:AInfty}.

Theorem \ref{theo:AInfty} below establishes that, under $L^{n/2}$ pinching conditions on the scalar curvature closely related to those in Theorem \ref{GHLn/2}, the volume densities $e^{nf}$ are \emph{uniformly $A_\infty$ weights} with respect to $g_0$.

 \begin{thm} \label{theo:AInfty}
Let  $(M^n,g_0)$ be a closed Riemannian manifold,  $R_0\in (0,\diam(M,g_0)]$ and $\delta, \Lambda>0$. Let $f : M\to \R$ be a smooth function such that the metric $g_f = e^{2f}g_0$ and its  associated Riemannian measure $d\mu_f=e^{nf}\mathrm{dvol}_{g_0}$ satisfy
\begin{enumerate}[i)]
\item $\int_{M} \left|\scal_{g_f}\right|^{\frac n2} \dv_{f}\le \Lambda$.
\item For any $x\in M$,  $\ds \left(\int_{\bB(x,R_0)} \left(\scal_{g_f}\right)_+^{\frac n2} \dv_{f}\right)^{\frac 2n}\le  \upalpha(n,2)-\delta.$
\end{enumerate}
Then $e^{nf}$ is an $A_\infty$-weight with respect to $g_0$, with constants only depending on $n, g_0, R_0, \delta$ and $\Lambda$. 
\end{thm}

The hypotheses as well as the conclusion of Theorem \ref{theo:AInfty} are scale invariant. We will see in Section \ref{sec:AInfty} that if we fix upper and lower volume bounds, having $A_\infty$ control for $e^{nf}$ implies $\cC^\alpha$ bounds for the distance associated to the metric $e^{2f}g$, for some $\alpha>1$. This implies as well uniform control on certain analytic quantities associated to the conformal metric. Denoting $g_f = e^{2f}g_0$, the fact that $e^{nf}$ is an $A_\infty$ weight implies that for any $g_0$-geodesic ball, its $g_f$-diameter is bounded from above in terms of its $g_f$-volume and that the measure $d\mu_f = e^{nf} d{\rm \mu}_{g_0}$ is doubling, with constants only depending on $\delta, \Lambda$ and $R_0$.

\medskip
 
Notice that in the second family of examples presented above, for which Theorem \ref{theo:AInfty} applies, the Riemannian volumes of the members of the family satisfy a \emph{uniform Strong $A_\infty$ bound} in the sense of \cite{DS}, see Definition \ref{def:sAInfW} in Section \ref{sec:SAInfty} below. This implies that for a sequence of metrics coming from these examples, the Gromov-Hausdorff limit $(M_\infty, d_\infty)$ is homeomorphic to $M$ and endows $M$ with a distance $d_{\infty}$ that is uniformly bi-H\"older to $(M,d_{g_0})$. In addition, the metrics in this family of examples have uniform Sobolev, Poincar\'e and isoperimetric inequalities. Our main result shows that this is always the case, provided the $L^{n/2}$-pinching of the scalar curvature is small enough at small scales.

\begin{thm}\label{theo:SAInfty}
Let  $(M^n,g_0)$ be a closed Riemannian manifold. Then there exists a constant $\Lambda_0 = \Lambda_0(g_0)$ with the following property: Let $g_f$ be given by $e^{2f} g_0$ with $f\in C^{2}(M)$. If there is a $R_0\in (0,\diam(M,g_0)]$ such that 
\begin{equation}
\quad \forall x\in M,  \  \int_{\bB(x,R_0)} \left|\scal_{g_f}\right|^{\frac n2} \dv_{f} \le  \Lambda_0
\end{equation} 
then  $e^{nf}$ is a strong $A_\infty$ weight with respect to $g_0$, with constants only depending on $g_0, R_0$ and $\Lambda_0$. Moreover, for all $\alpha\in (0,1)$, the distances $d_f$ and $d_0$ are uniformly $\alpha$-bi-H\"older with constants depending only on $g_0, R_0$, $\Lambda_0$ and $\alpha$.
\end{thm}

Before we continue, let us introduce some notation. Let $(M^n,g_0)$ be a closed Riemannian manifold, for $v, V, \delta, \Lambda>0$ and $R_0\in (0,\diam(M,g_0)]$, define $\cM_{v,V,R_0}^{\Lambda}$ as the set of smooth functions $f$ on $M$ such that the associated metric $g_f = e^{2f}g_0$ satisfies
\begin{enumerate}[i)]
\item $v\le \vol (M,g_f)\le V$
\item for any $x\in M\colon$ $\left(\int_{\bB(x,R_0)} \abs{\scal_{g_f}}^{\frac n2} \dv_{f}\right)^{\frac 2n}\le  \Lambda$.
\end{enumerate}

As mentioned previously, strong $A_\infty$-weight control coupled with volume bounds has many interesting geometric consequences. Some of them are listed in the following corollary.

\begin{cor}\label{coro:biHolder}
 Let  $(M^n,g_0)$ be a closed Riemannian manifold. Then there exists $\Lambda_0 = \Lambda(g_0)$ such that for all $R_0\in (0,\diam(M,g_0)]$ and any $v, V$ with $0<v\leq V$,  the set of metric spaces
$$\{(M, d_{g_f}) \; ; \; f\in \cM_{v,V,R_0}^{\Lambda_0}\}$$
is precompact for the Gromov-Hausdorff topology. In addition, any sequence in this set has a subsequence whose limit $(M, d_\infty)$ is $\alpha$-bi-H\"older to $(M, d_{g_0})$ for all $ \alpha\in (0,1)$. 
\end{cor}

It is worth noticing that the conclusion of Theorem \ref{theo:AInfty} is weaker than the conclusion of Corollary \ref{coro:biHolder} but holds under an explicit (and optimal, due to examples of \cite{CGW}) $L^{n/2}$-pinching of the scalar curvature. The constant $\Lambda(g_0, R_0)$ can be explicited in terms of the Ricci curvature of $g_0$ and its diameter. We could not produce a sequence of examples satisfying the hypotheses of Theorem \ref{theo:AInfty} which would not be uniformly Strong $A_\infty$.

\begin{rem}
Our proof of Theorem \ref{theo:SAInfty} shows that, under the same hypotheses, $e^{nf}$ is actually a \emph{stronger $A_\infty$ weight} in the sense of \cite{Semmes1}. This implies by Theorem 5.2 of \cite{Semmes1} that there exists $N>0$ such that any such $(M, g_f = e^{2f}g_0)$ can be embedded in $\R^N$ via a bi-Lipschitz embedding, with uniform bi-Lipschitz constant. We will not further discuss this fact here.
\end{rem}

\medskip

Eventhough our initial objective was to study conformal deformations of closed manifolds, our techniques also provide an interesting criterion to get strong $A_\infty$ on $(\mathbb R^n, \eucl)$.

\begin{thm}\label{theo:SAInfty-Rn}
Let $g = e^{2f} \eucl$ be a conformal deformation of the Euclidean metric on $\R^n$ such that:
\begin{itemize}
\item $\vol(\R^n, g) = +\infty$,
\item $\displaystyle \int_{\R^n} | \scal_g |^{n/2}d\mu_g < +\infty$.
\end{itemize}
Then $e^{nf}$ is a strong $A_\infty$ weight on $\R^n$ with respect to $\eucl$.
\end{thm}

As mentionned above, being a strong $A_\infty$ weight has many implications in terms of geometric control.
In particular, this theorem has the following striking geometric corollary, which seems to have remain unknown.

\begin{cor}\label{cor:Isoper-Rn}
Let $g = e^{2f} \eucl$ be a conformal deformation of the Euclidean metric on $\R^n$ such that:
\begin{itemize}
\item $\vol(\R^n, g) = +\infty$,
\item $\displaystyle \int_{\R^n} | \scal_g |^{n/2}d\mu_g < +\infty$.
\end{itemize}
Then there are positive constants $\uptheta, \upgamma$ such that any $g$-geodesic ball $B_g(x,r)$ satisfies
$$\uptheta^{-1} r^n\le \mu_g(B_g(x,r))\le \uptheta r^n.$$ 
Moreover, $(\R^n, g)$ satisfies the Euclidean isoperimetric inequality: for any bounded smooth open subset $\Omega\subset \R^n$ we have:
$$\upgamma \ \mu_g(\Omega)^{1-\frac 1n}\le \mu_g(\partial \Omega).$$
\end{cor}

\medskip

Let us conclude by an intricate question. Let $(g_k = e^{2f_k} g_0)_{k\in \N}$ be a sequence of conformal deformations of $g_0$ with $f_k\in \cM_{v,V,R_0}^{\Lambda_0}$, where $\Lambda_0$ satisfies the conclusion of Corollary \ref{coro:biHolder}. Then, up to extracting a subsequence, the sequence of associated metric spaces $(M, d_{g_k})$ converges in Gromov-Hausdorff topology to some distance $d_\infty$ on $M$, which is bi-H\"older equivalent to $d_0$.  Moreover, a subsequence of the Riemannian volume measures $\mu_k = d{\mu}_{g_k} = e^{n f_k}d{\mu}_{g_0}$ converges in the weak$-*$ topology to a positive measure $\mu_\infty$ on $M$. Since all the $\mu_k$ are uniformly strong $A_\infty$-weights with respects to $d{\rm vol}_{g_0}$, so is $\mu_\infty$: there exists $e^{n f_\infty}\in L^1(M)$ such that $\mu_\infty= e^{n f_\infty}d{\rm vol}_{g_0}$ and $\mu_\infty$ is still a strong $A_\infty$-weight with respect to $d{\rm vol}_{g_0}$. Hence the weight $e^{n f_\infty}$ induces another distance $d_{f_\infty}$ on $M$, analogous to a conformal Riemannian distance (see Section \ref{sec:LimDist} below or \cite{Semmes1} for a precise definition).

A natural question is: \emph{Under which assumption do the distances $d_\infty$ and $d_{f_\infty}$ coincide?} We could not give a full answer to this question. We know that they do coincide in  Example 2 presented above. On the other hand, in the last section of this paper we will describe another example which shows that uniform strong $A_\infty$-weight control is not enough to establish this equality in general. Nevertheless, volume bounds and a $L^{n/2}$ bound on the scalar curvature give a much stronger control on the geometry than just strong $A_\infty$-weight controls. Thus, the distances could still coincide under the integral pinching of the scalar curvature under consideration.

\medskip

From now on, given a background Riemannian manifold $(M, g_0)$, for any conformal deformation $g_f = e^{2 f} g_0$, we will denote by $d_f = d_{g_f}$ the associated Riemannian distance and by $d\mu_f = d{\rm vol}_{g_f}$ the associated Riemannian volume measure. In particular, $d_0$ and $\mu_0$ are the distances and volume measure associated to $g_0$, respectively.

\medskip

{\bf Acknowledgements:} The authors are grateful for the support of the ANR grants ACG (ANR-10-BLAN 0105) and CCEM (ANR-17-CE40-0034) which made this collaboration possible. C.L. Aldana was also supported by the Fonds National de la Recherche, Luxembourg 7926179.

\section{Integral bounds for the scalar curvature and convergence of distances} 

This first section is devoted to the proof of Theorem \ref{GHLn/2}. This theorem is an immediate consequence of Propositions \ref{GHPrecompactness} and \ref{Integimprove} given below. The proof will show that for any sequence of metrics $(g_n)_{n\in \N}$ in $\cM_{V}^{\delta}$, the sequence of distances $(d_{g_n})_{n\in \N}$ converges in some H\"older topology up to extracting a subsequence. This implies the precompactness of the sequence of metric spaces $(M, d_{g_n})_{n\in \N}$ for the Gromov-Hausdorff topology, as shown by the following lemma which will be used several times in the sequel.

\begin{lem}\label{lem:Holder-Implies-GH}
Let $(X,d_0)$ be a compact metric space, $\alpha\in (0,1)$, $C>0$, and $(d_n)_{n\in \mathbb N}$ be a sequence of distances on $M$ which is bounded in the $\alpha$-H\"older topology with respect to $d_0$, i.e. such that for all $x, x',y, y'\in M$ and all $n\in \mathbb N$,
$$
|d_n(x,y) - d_n(x',y')| \leq C\left(d_0(x,x')^\alpha + d_0(y,y')^\alpha\right).
$$
Then up to extracting a subsequence, $d_n$ is converging uniformly and in the $\alpha'$-H\"older topology for all $\alpha'\in (0, \alpha)$ to a map $d_\infty : M\times M\to [0, +\infty)$. Moreover, the sequence of metric spaces $(M, d_n)$ converges in the Gromov-Hausdorff topology to the metric space $(M_\infty, d_\infty)$ where $M_\infty = M/\sim$, with $x\sim y$ if and only if $d_\infty(x,y) = 0$. 
\end{lem}

\begin{proof}
Let $(d_n)_{n\in \mathbb N}$ be a sequence of distances on $M$ satisfying the above hypotheses. By Arzela-Ascoli Theorem, it has a subsequence which converges in the $\cC^0$-uniform topology on $M\times M$ to some non-negative map $d_\infty:M\times M\to [0, +\infty)$, satisfying the triangular inequality. By construction, $d_\infty$ induces a distance on $M_\infty = M/\sim$, with $x\sim y$ if and only if $d_\infty(x,y) = 0$. The $\cC^0$-uniform convergence of the $d_k$ implies that the sequence of metric spaces $(M, d_n)$ converges in the Gromov-Hausdorff topology to the metric space $(M_\infty, d_\infty)$, see for instance \cite{BBI} p.260, Example 7.4.4.

An elementary computation shows that the convergence of $(d_n)$ is actually uniform in all $\cC^{\alpha'}$ H\"older spaces for all $0 < \alpha'< \alpha$, which concludes the proof of our lemma.
\end{proof}

\begin{prop}\label{GHPrecompactness}
Let $(M,g_0)$ be a compact Riemannian manifold and let $(g_k = e^{2f_k}g_0)_{k\geq 0}$ be a sequence of smooth conformal metrics. Assume that there exist $p>n$ and $A>0$ such that for all $k\in \N$,
$$\int_M e^{p f_k} d{\mu}_0 \leq A.$$
Then the sequence of distances $\ds \left( d_{f_k} : M\times M\to [0, +\infty)\right)_{k\in \N}$ is precompact in the $\alpha$-H\"older topology for all $\alpha \in (0, 1 - \frac n p)$. In particular, the sequence of metric spaces $(M, d_{f_k})_{k\in \N}$ is precompact for the Gromov-Hausdorff topology.
\end{prop}

\begin{proof}
Let $x\in M$ be fixed. Recall that for any smooth Riemannian metric $g$ on $M$, $\norm{\nabla^gd_g(x,.)}_g = 1$ almost everywhere, see e.g. \cite[Prop. 4.8]{Sakai}. Let $g_f = e^{2f} g_0$ be conformal metric, we have hence for almost all $y\neq x$,
$$
\norm{\nabla^{g_0}d_f(x,.)}_{g_0}(y) = e^{f(y)}.
$$
Therefore, there exists $C>0$ such that
\begin{equation}\label{W1pBound}
\int_M e^{p f} d{\mu}_{0} \leq A \quad \Rightarrow \quad \norm{d_{f}(x,.)}_{W^{1,p}(M, g_0)} \leq CA,
\end{equation}
where $W^{1,p}(M, g_0)$ is the Sobolev space of functions on $M$ with distributional derivative in $L^p(g_0)$.

Now, let $(g_k = e^{2f_k}g_0)_{k\geq 0}$ be a sequence of smooth conformal metrics. Assume that there exist $p>1$ and $A>0$ such that for all $k\in \N$,
$$\int_M e^{p f_k} d{\mu}_{0} \leq A.$$
It follows from (\ref{W1pBound}) that the sequence of distances $(d_{f_k}(x,.))_{k\in \N}$ is bounded in $W^{1,p}(M, g_0)$. Therefore, since for all $k\in \N$, $d_{f_k}(x,.)$ vanishes at $x$, by Sobolev embeddings the sequence is bounded in $\cC^\alpha$ with $\alpha = 1- \frac n p$ by some constant $C'$ independent of $x$. 

We have hence, for all $x, x',y, y'\in M$ and all $n\in \mathbb N$,
\begin{eqnarray*}
|d_{f_k}(x,y) - d_{f_k}(x',y')| & = & |d_{f_k}(x,y) - d_{f_k}(x',y) + d_{f_k}(x',y) - d_{f_k}(x',y')|\\
& \leq & C'\left(d_0(x,x')^\alpha + d_0(y,y')^\alpha\right).
\end{eqnarray*}

By Lemma \ref{lem:Holder-Implies-GH}, this implies Proposition \ref{GHPrecompactness}.

\end{proof}

We now show that any conformal deformation of $g_0$ whose scalar curvature has a positive part with $L^{n/2}$ norm strictly less than the round sphere has a conformal factor that satisfies the integrability condition required in Proposition \ref{GHPrecompactness}.

\begin{prop}\label{Integimprove} 
Let $(M^n, g_0)$ be a closed Riemannian manifold. Let $V,\delta$ be posi\-tive real numbers. There exists $ A = A(V, \delta)>0$ such that for all $\e>0$ small enough, if  
\begin{enumerate}[i)]
\item $\vol(M,g_f)\le V,$
\item$ \norm{(\scal_{g_f})_+}_{L^{\frac{n}{2}}(M, \dv_{f})}\le\upalpha(n,2)-\delta$
\end{enumerate} 
then 
\begin{equation}\label{eq:IntegSup}
\int_M e^{n(1+\upepsilon) f} \dv_{g_0}\le A.
\end{equation}
\end{prop}

\begin{proof}

We are going to use the following optimal Sobolev inequality du to E. Hebey and M. Vaugon (\cite{Hebey-Vaugon}) (this inequality had been conjectured by T. Aubin \cite{Aubin}):
if we write again $\upalpha(n,2)=\sigma_n^{\frac 2n} n(n-1)$, then there is a constant $B$ such that for all $\varphi\in \cC^\infty(M)$, 

\begin{equation} \label{eq:HabeyVaugon}
\upalpha(n,2)\left(\int_M \varphi^{\frac{2n}{n-2}}{\rm d}\mu_0\right)^{1-\frac2n}\le \frac{4(n-1)}{n-2}\int_M |d \varphi|_{g_0}^2\dv_{g_0} +\int_MB\varphi^2\\dv_{g_0} .
\end{equation}
If $g_f=u^{\frac{4}{n-2}} g_0$ if follows from \cite{Yam60} that
\begin{equation}\label{Yamabe}
\scal_{g_f} u^{\frac{4}{n-2}}\,u=\frac{4(n-1)}{n-2} \Delta_{g_0} u+ \scal_{g_0} u.
\end{equation} 

In particular for any $\upepsilon>0$ we get 
\begin{multline*}
\frac{4(n-1)}{n-2}\Delta_{g_0} u^{1+\upepsilon}+ (1+\upepsilon)\scal_{g_0} u^{1+\upepsilon}+(n-1)\frac{4\upepsilon(1+\upepsilon)}{n-2}u^{\upepsilon-1}|du|^2_{g_0}\\=(1+\upepsilon)\scal_{g_f} u^{\frac{4}{n-2}}\,u^{(1+\upepsilon)}.
\end{multline*}
Using the Hebey and Vaugon Sobolev inequality for $\varphi=u^{(1+\upepsilon)}$, we get
\begin{eqnarray*}
 & & \upalpha(n,2)\left(\int_M u^{(1+\upepsilon)\frac{2n}{n-2}}\dv_{g_0}\right)^{1-\frac2n}\\
 & & \hspace{1cm} \le \frac{4(n-1)}{n-2}\int_M u^{(1+\upepsilon)}\Delta_{g_0} u^{1+\upepsilon}\dv_{g_0}
 +B\int_M u^{2(1+\upepsilon)} \dv_{g_0}\\
& & \hspace{1cm}  \le (1+\upepsilon)\int_M(\scal_{g_f})_+ u^{\frac{4}{n-2}}\,u^{2(1+\upepsilon)}\dv_{g_0} \\
& & \hspace{1.5cm} +((1+\upepsilon)S_0+B)\int_Mu^{2(1+\upepsilon)}\dv_{g_0},
\end{eqnarray*}

where $S_0=\|\scal_{g_0}\|_{L^\infty}.$
When $\upepsilon\le 2/(n-2)$, using the H\"older inequality we get
\begin{align*}&\left((\upalpha(n,2)-(1+\upepsilon)\|(\scal_{g_f)_+}\|_{L^{\frac{n}{2}}(M, \dv_{f})} \right)\left(\int_M u^{(1+\upepsilon)\frac{2n}{n-2}}\dv_{g_0}\right)^{1-\frac2n}\le \\
&\hspace{1cm}((1+\upepsilon)S_0+B)\left(\vol_{g_0}(M)\right)^{\frac{2-(n-2)\upepsilon}{n}}\left(\vol_{g_f}(M)\right)^{(1+\upepsilon)\frac{n-2}{n}},\end{align*}

which gives with $e^{2f} = u^{\frac{4}{n-2}}$,
$$
\int_M e^{n(1+\upepsilon) f} \dv_{g_0} \leq \frac{((1+\upepsilon)S_0+B)\left(\vol_{g_0}(M)\right)^{\frac{2-(n-2)\upepsilon}{n}}}{\delta - \upepsilon(\upalpha(n,2)-\delta)} V^{(1+\upepsilon)\frac{n-2}2},
$$
which is precisely (\ref{eq:IntegSup}) when $\upepsilon$ is small enough, say $\upepsilon < \frac{\delta}{\upalpha(n,2)-\delta}$.

\end{proof}

 Notice that under our very weak hypotheses, the limit space $(M_\infty, d_\infty)$ obtained as the Gromov-Hausdorff limit of the sequence $(M, d_{g_k})_{k\in\N}$ may be very singular; the singularities depending on the set where $d_\infty$ vanishes. Even if the volume of $(M, g_k)$ remains constant (recall that the $L^{n/2}$ bound of the scalar curvature is invariant under rescaling), we cannot a priori prevent collapsing of some open part of the manifold. In the sequel, we look for stronger conditions on the scalar curvatures and the volumes so that we get a better geometric picture of the limit.

Let us first note that if we add to the hypotheses in Theorem \ref{GHLn/2} a lower bound on the volume and a uniform $L^p$ bound on the \emph{scalar curvature} (and not the full Riemann tensor) with $p>\frac n 2$, then we get $\cC^{\alpha}$-precompactness for the conformal factors:

\begin{cor}\label{Gurskyimprove} Let $(M^n,g_0)$ be a closed Riemannian manifold. Let $v,V,\Lambda,\delta$ be positive real numbers and let $p>n/2$. 
Then the set 
$\cM_{v,V}^{p,\Lambda,\delta}$ of smooth functions $f\in \cC^\infty(M)$ such that $g_f=e^{2f}g_0$ satisfies
\begin{enumerate}[i)]
\item $v\le \vol(M,g_f)\le V,$
\item $ \int_M\left|\scal_{g_f}\right|^p \dv_{f}\le \Lambda$
\item$ \norm{(\scal_{g_f})_+}_{L^{\frac{n}{2}}(M, \dv_{f})}\le\sigma_n^{\frac 2n} n(n-2)-\delta$
\end{enumerate} is precompact in $\cC^{\alpha}(M)$ for $0<\alpha<\frac{2p-n}{p}$.
\end{cor}

\begin{proof} 
It follows from Proposition \ref{Integimprove} that for all metrics $g_f = e^{2f}g_0$ in $\cM_{v,V}^{p,\Lambda,\delta}$, the volume density $e^{nf}$ satisfies a $L^{1+\epsilon}$ uniform bound given by (\ref{eq:IntegSup}). 

Combining this uniform stronger integrability of $e^f$ together with $L^p$ bounds for the scalar curvature and volume bounds, a Harnack inequality due to Trudinger in \cite{Trudinger} together with elliptic regularity will give $\cC^{\alpha}$-precompactness for the space of metrics for $0<\alpha<\frac{2p-n}{p}$. This is explained with details in Section 2.2 of \cite{HM}.
\end{proof}

In particular, any sequence of metrics satisfying some $L^p$ bound on the scalar curvature, $p>n/2$, and the above $L^{n/2}$ bounds on its positive part, converges (up to extraction of a subsequence) in some $\cC^\alpha$-topology to a Riemannian metric $g_\infty = e^{2f_\infty} g_0$, where $f_\infty\in \cC^{\alpha}(M)$. This slightly improves the well known result by M. Gursky given in \cite{Gursky}.

\medskip
The rest of our paper is devoted to the study of sequences of conformal deformations with fixed volume bounds and a mere $L^{n/2}$ bound on the scalar curvature. It follows from \cite{CGW} that in this setting, one cannot get $\cC^0$ precompactness for the sequence of conformal factors (or, equivalently, for the sequence of Riemannian metrics). Nevertheless, we will show that the sequence of associated measured metric spaces and the limit metric spaces $(M_\infty, d_\infty)$ thus obtained satisfy many uniform geometrical properties.

\section{Optimal $L^{n/2}$ bounds for the scalar curvature and $A_\infty$ weights}\label{sec:AInfty}

In this section, we study sequences $(g_k = e^{2f_k} g_0)_{k\in \N}$ of conformal metrics with a uniform upper bound on the volumes, such that the $L^{n/2}$-norms of the positive part of the scalar curvature on small spheres is uniformly bounded by the one of the standard sphere, and which moreover have some uniform $L^{n/2}$ bound for the total scalar curvature. We will show that the Riemannian volume densities $\{ e^{n f_k}\}_{k\in \N}$ are then \emph{uniformly $A_\infty$-weights} with respect to $g_0$, see Definition \ref{defi:AInfty}.

To be uniformly $A_\infty$-weights has several nice geometrical consequences. First of all, the sequence of volume measures $(d\mu_{f_k})_{k\in \N}$ is uniformly doubling with respect to the initial distance $d_{0}$. In particular, this rules out the kind of blowups presented in Example 1 in our introduction. More precisely, we will see that the $A_\infty$ control implies that the Riemannian distances $d_{f_k}$ are uniformly controlled in terms of the volumes $d\mu_{f_k}$, locally. We will show as well Gromov-Hausdorff precompactness for the set of conformal metrics satisfying a uniform $A_\infty$-weight bound, independently of any control on the curvature.

\medskip

In the afore mentioned work of Gursky \cite{Gursky}, a key step to get $\cC^{\alpha}$ compactness for a sequence of conformal metrics $(g_k = e^{2f_k} g_0)_{k\in \N}$ whose Riemann curvatures satisfy a uniform $L^p$ bound ($p>n/2$) is to show that if $g_0$ is not conformal to the standard sphere, then the conformal factors in such sequences satisfy a uniform Harnack inequality: there exists $C>1$ such that for all $k\in \N$,
$$\frac{\sup e^{f_k}}{\inf e^{f_k}} \leq C.$$

Uniform volume bounds and (even very small) $L^{n/2}$ bounds on the scalar curvature cannot ensure such Harnack inequality, as shows Example 2 in our introduction. A uniform $A_\infty$ control on the weight $e^{n f_k}$ is actually a local $L^q$ integral version of this Harnack inequality for some $q>1$.

\subsection{$A_\infty$ weights and conformal metrics}

We now give a presentation of $A_\infty$ weights on a closed Riemannian manifold. The reader can find a good exposition and more details in \cite[chapter V]{Stein}. Then, we will draw some geometric consequences for conformal metrics having such $A_\infty$ control.

Let $(M^n,g_0)$ be a closed Riemannian manifold. Throughout this document we will denote by
$$\bB(x,r) :  \text{ the } g_0\text{-geodesic ball centered at } x \text{ and of radius } r. $$ 
We will sometimes write $\bB$ for a $g_0$-geodesic ball whose radius will be denoted by $r(\bB)$. For $\uptheta>0$, $\uptheta \bB$ will be the ball with the same center and with radius $\uptheta r(\bB)$. $D$ will denote the diameter of $(M^n,g_0)$. The measure $\mu_0$ will still be the $g_0$ Riemannian volume. If $E\subset M$ is a measurable set and $f\in L^1(E,d\mu_0)$ the average of $f$ over $E$ will be denoted by $$\fint_E fd\mu_0=\frac{1}{\mu_0(E)}\int_E fd\mu_0.$$

\begin{defi} \label{defi:AInfty}
We say that a non negative function $w\in L^1(M,d\mu_0)$ is a \emph{$A_\infty$ weight with respect to $g_0$} if one of the following equivalent properties is satisfied:
\begin{enumerate}[i)]
\item There is a $q>1$ and a constant $C$ such that for any geodesic ball $\bB\subset M$, the following reverse H\"older inequality with exponent $q$ holds:
$$\left(\fint_{\bB} w^qd\mu_0\right)^{\frac{1}{ q}}\le C \fint_{\bB} wd\mu_0.$$

\item There is $p>1$ and a constant $C$ such that for any geodesic ball $\bB\subset M$:
$$\left(\fint_{\bB} wd\mu_0\right)\left(\fint_{\bB} w^{-\frac{1}{p-1}}d\mu_0\right)^{p-1}\le C.$$
\emph{If $w$ satisfies this condition, it is called $A_p$-weight.}
\item There are constant $\delta,\epsilon\in (0,1)$ such that for any geodesic ball $\bB\subset M$  and any $E\subset M$ if 
$E\subset \bB$ satisfies $\mu_0(E)\le \delta \mu_0(\bB)$ then $\int_E wd\mu_0\le \epsilon \int_{\bB} wd\mu_0.$
\item There are constants $\alpha>1$ and $C>0$ such that for any geodesic ball $\bB\subset M$  and any $E\subset \bB$:
$$\frac{1}{C}\left(\frac{\mu_0(E)}{\mu_0(\bB)}\right)^\alpha\le \frac{\int_E wd\mu_0}{\int_{\bB} wd\mu_0}\le C\left(\frac{\mu_0(E)}{\mu_0(\bB)}\right)^{\frac 1\alpha.}$$
\end{enumerate}
\end{defi}

\begin{rem} \label{Rmks:Sb-to-Bb}
The different constants in the different definitions are mutually controlled.
In particular, we will use the following fact:  If $w$ is an $A_p$ weight satisfying part \emph{ii)} above for all $p>1$, then it fulfills reverse H\"older inequality \emph{i)} with exponent $q$ for all $q>1$, with mutually controlled constants. See \cite[p.196-203]{Stein} for a proof. 
\end{rem}

It is possible to show that the properties in Definition \ref{defi:AInfty} hold for all geodesic balls if and only if they hold for geodesic balls of radius less than some fixed $R_0>0$.  We will prove this statement for the Reverse H\"older inequality. In order to do that, we need the following general result.
 
 \begin{lem}\label{lem-doublingrR} If $(X,d,\mu)$ is a metric measure space whose balls $B$ with radius less than $R_0$ satisfy the doubling condition:
 $$\mu(2B)\le \uptheta \mu(B)$$ then for any $R_1\ge R_0$ there is a constant $\uptheta'$ depending only on $\uptheta$ and $R_1/R_0$ such that 
 any ball $B$ with radius less than $R_1$ satisfies the doubling condition:
 $$\mu(2B)\le \uptheta' \mu(B)$$ 
 \end{lem}
 \proof Indeed let $B(x,R)$ a ball of radius $R\in [R_0,R_1]$, then from \cite[Lemma 3.10]{Carron2016}, we have 
$$\mu(B(x,2R))\le \uptheta^{50+50\frac{2R}{R_0}}\mu_f(B(x,R_0))\le \uptheta^{50+50\frac{2R_1}{R_0}}\mu(B(x,R)).$$ 
 \endproof
 
 \begin{lem}\label{lem:Ap at all scales}
Let $q>1$, and $C, R_0>0$ be fixed. Let $w$ be a uniformly $(q,R_0,C)$-$A_\infty$ weight, i.e. every
geodesic ball $\bB$ of radius $r\le R_0$ satifies:
$$
\left(\fint_{\bB} w^qd\mu_0\right)^{\frac{1}{q}}\le C \fint_{\bB} wd\mu_0.
$$
Then the measure $w d\mu_0$ has the doubling volume property and there exists $C' = C'(q,C,R_0, g_0)>0$ such that for all balls $\bB$ of radius $r\geq 0$, one has
\begin{equation}
\left(\fint_{\bB} w^qd\mu_0\right)^{\frac{1}{q}}\le C' \fint_{\bB} wd\mu_0.
\end{equation}
\end{lem}

\begin{proof}
Assume that $w$ satisfies the hypotheses of the lemma for some $q>1$ and $C,R_0>0$. We first show that $w d\mu_0$ has the doubling volume property. Let us first consider geodesic balls of radius less than $R_0>0$. We start by showing that there is $\uptheta>0$ such that for the measure $d\omega=wd\mu$ we have the doubling property:

$$
r(\bB)\le R_0/2\Longrightarrow \omega(2\bB)\le \uptheta \omega(\bB).
$$
Indeed using both H\"older and Reverse H\"older inequalities we easily get for any $\tau\in (0,1):$
$$\frac{\omega(\bB\setminus \tau\bB)}{\omega(\bB)}\le C\left(\frac{\mu_0(\bB\setminus \tau\bB)}{\mu_0(\bB)}\right)^{1-\frac 1q}.$$
The manifold $M$ is compact hence there is a constant $B$ such that for any geodesic ball :
$$\mu_0(\bB\setminus \tau\bB)\le B (1-\tau) \mu_0(\bB).$$
We choose $\tau\in (0,1)$ such that 
$C\left(B (1-\tau)\right)^{1-\frac1 q}=\frac 12$ and we get $\omega(\bB)\le 2  \omega(\tau\bB)$,
therefore
$\omega(\bB)\le 2^N  \omega(\tau^N\bB)$.
Choosing now $N$ such that $\ds \tau^N\le \frac 12< \tau^{N-1}$ gives the result with $\uptheta=2^N .$ 

As the diameter of $(M,g_0)$ is finite, we deduce  from the \lref{lem-doublingrR}  that the measure $\omega$ is doubling: there is a constant $\uptheta$ such that for any ball $\bB\subset M$:
$$\omega(2\bB)\le \uptheta \omega(\bB).$$

We can now show that the reverse H\"older inequality holds for any ball.
Let $\bB$ be a geodesic ball of radius $r\in (R_0,D]$, then we can find a minimal family of geodesic balls $\bB_\alpha =\bB(x_\alpha,R_0/2)$ of radius $R_0/2$ such that $\bB\subset \cup_\alpha\bB_\alpha$, and the balls $\bB(x_\alpha,R_0/4)$ are disjoint, hence include in $2\bB$. Since $w$ satisfies reverse H\"older on the balls $\bB(x_\alpha,R_0/2)$, we get
$$
\left(\int_{\bB_\alpha} w^qd\mu_0 \right)^{\frac 1q} \le \frac{C\mu_0(\bB_\alpha)^{\frac1q}}{\mu_0(\bB_\alpha)}\int_{\bB_\alpha} w d\mu_0 = C\mu_0(\bB_\alpha)^{\frac1q-1}\int_{\bB_\alpha} w d\mu_0.
$$
From the doubling condition on $(M,d_0,\mu_0)$ there is a constant such that 
$$\mu_0(\bB)\le C\mu_0(\bB_\alpha).$$
Hence
$$
\left(\int_{\bB_\alpha} w^qd\mu_0 \right)^{\frac 1q} \le C\mu_0(\bB)^{\frac1q-1}\int_{\bB_\alpha} w d\mu_0 $$

\begin{align*}
\int_{\bB} w^qd\mu_0
&\le \sum_{\alpha} \int_{\bB_\alpha} w^qd\mu_0\\
&\le C \sum_{\alpha} \mu_0(\bB_\alpha)^{1-q}\left(\int_{\bB_\alpha} wd\mu_0\right)^q\\
&\le  C\mu_0(\bB)^{1-q}\sum_{\alpha}\left(\int_{\bB_\alpha} w d\mu_0\right)^q \\
&\le C\mu_0(\bB)^{1-q} \uptheta^q  \sum_{\alpha}\left( \omega(\frac12\bB_\alpha)\right)^{q}\\
&\le C\mu_0(\bB)^{1-q} \uptheta^q \left(\sum_{\alpha}\omega(\frac12\bB_\alpha)\right)^q\\
&\le C\mu_0(\bB)^{1-q} \uptheta^q  \omega(2\bB)^q  \le C\mu_0(\bB)^{1-q} \uptheta^{2q}  \omega(\bB)^q,
\end{align*}
from which the conclusion follows.
\end{proof}

Let now $f\in \cC^{\infty}(M)$ and $g_f = e^{2 f}g_0$ be a Riemannian metric conformal to $g_0$. A key geometric consequence of having $A_\infty$ control on the weight $e^{n f}$ is given in the following lemma.

\begin{lem}\label{lowerd} Let $f : M\to \R$ be a smooth map. 
Assume that $w=e^{nf}$ is an $A_\infty$ weight  satisfying the above reverse H\"older inequality with constants $C>0$ and $q>1$. Then there is a constant $B$, depending only on $C, q$ and $(M,g_0)$, such that for any $x,y\in M$  :
$$d_{f}(x,y)^n\le B\int_{\bB(x,d_0(x,y))} e^{nf}\mathrm{d}\mu_0.$$
\end{lem}

This lemma has been shown on the Euclidean space by G. David and S. Semmes, see \cite[inequality (1.2)]{DS} and  \cite[Proposition 3.12b]{Semmes1}.

\begin{proof} By Sobolev inequality (see for instance \cite{GT} p.148), there is a constant $C>0$ depending only on the geometry of $(M,g_0)$ and $q$, such that for $p=qn$, we have for any balls $\bB\subset M$ and any $\varphi\in W^{1,p}(2\bB)$ and $x,y\in \bB$:
$$|\varphi(x)-\varphi(y)|\le C r(\bB)^{1-\frac{n}{p}}\left(\int_{2\bB} |d\varphi |^pd\mu_0\right)^{\frac 1p}.$$
If we use this for the function $\varphi=d_f(x,.)$ and $\bB$ a ball centered at some point $m$ with $d_0(x,m)=d_0(m,y)=d_0(x,y)/2$ and with radius $\frac 34 d_0(x,y)$. we will get 
$$d_f(x,y)\le C r(\bB)^{1-\frac{n}{p}}\left(\int_{2\bB} e^{qnf}d\mu_0\right)^{\frac{1}{qn}}$$
Then the conclusion follows from the reverse H\"older inequality and the fact that the $\mu_0$ measure of $2\bB$ is comparable to $r(\bB)^n$.
\end{proof}

Let us show now that the space of conformal metrics with uniform $A_\infty$ bounds is precompact in the Gromov-Hausdorff topology. We will not need this result in the sequel since Theorem \ref{GHLn/2} will already ensure Gromov-Hausdorff precompactness once we fix a suitable $L^{n/2}$ bound on the scalar curvature. Nevertheless, this is an easy consequence of the measure being $A_\infty$ weights which seems to have remained unnoticed.

\begin{prop}\label{prop:AInfty-Compact}
Let $(M^n,g_0)$ be a closed Riemannian manifold let  $C,V>0$, $R_0\in (0,\diam(M,g_0)]$ and $q>1$.

Let $\cM_{V,R_0,q,C}$ be the set of smooth functions such that $f\in \cM_{V,R_0,q,C}$ satisfies the following conditions: 
\begin{enumerate}[i)]
\item for the Riemanian metric $g_f = e^{2f} g_0$, $\vol(M,g_f)\le V;$
\item $w=e^{nf}$ is $(R_0, q, C)$-uniformly $A_\infty$ i.e. for any geodesic ball $\bB$ of radius less or equal than  $R_0$ we have:
$$\left(\fint_{\bB} w^qd\mu_0\right)^{\frac{1}{q}}\le C \fint_{\bB} wd\mu_0.$$
\end{enumerate}
Then the set of distances $\{ d_{f}, f\in \cM_{V,R_0,q,C}\}$ is pre-compact in the $\cC^\alpha$ topology for all $\alpha\in (0, 1 - \frac 1 q)$. In particular, the set of metric spaces $\{ (M,d_{f}), f\in \cM_{V,R_0,q,C}\}$ is precompact for the Gromov-Hausdorff topology.
\end{prop}

\begin{proof}
By Lemma \ref{lem:Ap at all scales}, there exists $C'>0$ such that for all $f\in \cM_{v,V,R_0,q,C}$ and all geodesic balls $\bB$ (of any radius), 
$$\left(\fint_{\bB} w^qd\mu_0\right)^{\frac{1}{q}}\le C \fint_{\bB} wd\mu_0.$$
In particular, if the radius of $\bB$ is larger than the diameter of $(M,g_0)$, we get that for all $f\in \cM_{v,V,R_0,q,C}$,
$$
\int_M e^{nq f} \leq C' (\mbox{Vol}(M,g_0))^{\frac 1 q - 1} V.
$$
The conclusion follows then immediately from Proposition \ref{GHPrecompactness}.
\end{proof}

\subsection{$A_\infty$ control from $L^{n/2}$ bound on the scalar curvature}

We now show that a $L^{n/2}$ bound on the scalar curvature, with $L^{n/2}$ bound on the positive part of the scalar curvature uniformly smaller than what it is for the standard sphere at some (even small) fixed scale, implies a uniform $A_\infty$ control of the conformal factor. Once again, it follows from the Example 1 presented in our introduction that these integral bounds are optimal to get such $A_\infty$ control. Theorem \ref{theo:Conformal-AInfty} below implies Theorem \ref{theo:AInfty} of the introduction and explicit its result.

\begin{thm} \label{theo:Conformal-AInfty}
Let  $(M^n,g_0)$ be a closed Riemannian manifold, $\delta, \Lambda>0$ and $R_0>0$. 
Then there exist 
$p = p(n, g_0, R_0, \Lambda)>0$ and $C = C(n, g_0, R_0, \delta, \Lambda)>0$ such that for any smooth function $f : M\to \R$ whose associate metric $g_f = e^{2f}g_0$ satisfies
\begin{enumerate}[i)]
\item $\int_{M} \left|\scal_{g_f}\right|^{\frac n2}d\mu_f\le \Lambda$, and
\item for any $x\in M\colon$ $\ds \left(\int_{\bB(x,R_0)} \left(\scal_{g_f}\right)_+^{\frac n2}d\mu_f\right)^{\frac 2n}\le  \upalpha(n,2)-\delta,$
\end{enumerate}
we have that for any $g_0$-geodesic ball $\bB$ of radius less that $R_0/2$, 
$$\left(\int_\bB e^{nf}d\mu_0\right) \left(\int_\bB e^{-\frac{nf}{p}}d\mu_0 \right)^p \leq C.$$
In particular, $e^{nf}$ is an $A_\infty$-weight with respect to $g_0$, with constants only depending on $n, g_0, R_0, \delta$ and $\Lambda$.
\end{thm}

\proof The proof of this theorem was directly inspired by the proof of the Harnack inequality for positive solution of second order elliptic equation, see \cite[theorem 8.20]{GT} and its proof.

Let $R_0\in (0,\diam(M,g_0])$, $\delta, \Lambda>0$ and let $f : M\to\R$ be a smooth function satisfying the hypotheses above. We write again $g_f = u^{\frac{4}{n-2}}g_0$, i.e. $e^{f} = u^{\frac{2}{n-2}}$. The proof is done in two steps.

\begin{lem}\label{lem:AInfty1} Under the same hypothesis of Theorem \ref{theo:Conformal-AInfty}, there exist $\epsilon = \epsilon(n, g_0, R_0, \Lambda)>0$ and $C =  C_1(n, g_0, R_0, \Lambda)>0$ such that for any $g_0$-geodesic ball $\bB$ of radius less than $R_0/2$,
$$\left(\fint_{\bB} u^{\epsilon}\mathrm{d}\mu_0\right)\left(\fint_{\bB} u^{-\epsilon}\mathrm{d}\mu_0\right)\le C.$$
\end{lem}

\begin{proof}
Since we consider
$g_f=u^{\frac{4}{n-2}} g_0$, let us recall equation \ref{Yamabe}:
$$ \scal_{g_f} u^{\frac{4}{n-2}}\,u=\frac{4(n-1)}{n-2} \Delta_{g_0} u+ \scal_{g_0} u$$ 
Set $V=\scal_{g_f} u^{\frac{4}{n-2}}$, $W=V-\scal_{g_0} $ and $\lambda_0=\|\scal_{g_0}\|_{L^\infty}.$ We have
$$\frac{4(n-1)}{n-2}\Delta_0 \ln(u)=W+ \frac{4(n-1)}{n-2} |d\ln(u)|^2_{g_0}.$$
Now let $\bB$ and $\xi$ be a Lipschitz cut-off function with support in $2\bB$ such that
$\xi=1$ on $\bB$ and $|d\xi|\le r(\bB)^{-1}.$
Then 
$$
\int_M \xi^2\left(W+ \frac{4(n-1)}{n-2} |d\ln(u)|^2_{g_0}\right) \mathrm{d}\mu_0 =\frac{4(n-1)}{n-2}\int_M 2\xi \langle d\xi, d\ln(u)\rangle \mathrm{d}\mu_0
$$
$$
\le 4\int_M \xi^2|d\ln(u)|^2_{g_0}\mathrm{d}\mu_0+ \frac{16(n-1)^2}{(n-2)^2}\int_M |d\xi|^2 \mathrm{d}\mu_0
$$

By H\"older inequality, we have
$$\int_M \xi^2 (-V){\rm d}\mu_0 \leq \norm{Scal_{g_f}^-}_{L^{n/2}(2\bB)}\mu_0(2\bB)^{\frac{n-2}{2}},$$
Moreover, for $g_0$ there is an uniform constant $\omega$ such that for any $r\in (0,diam(M,g_0)]$ and any $x\in M$:
\begin{equation} \label{eq:rel-vol-rad-mu0}
\frac 1\omega r^n\le \mu_0(\bB(x,r))\le \omega r^n.
\end{equation} 

Hence we get
$$
\left(\frac{4(n-1)}{n-2}-4\right)\int_{\bB} |d\ln(u)|^2_{g_0}\mathrm{d}\mu_0\le C\lambda_0r(\bB)^n+C\left(1+\norm{Scal_{g_f}^-}_{L^{\frac n2}(2\bB)}\right)\mu_0(2\bB) r(\bB)^{-2}.
$$
In particular we have
\begin{equation}\label{BMO}
r(\bB)^2\fint_{\bB} |d\ln(u)|^2_{g_0}\mathrm{d}\mu_0\le C\left(1+\norm{Scal_{g_f}^-}_{L^{\frac n2}(2\bB)}\right).\end{equation}
There is a $r_0>0$ such that all the geodesic balls $\bB(x,r_0)$ are almost Euclidean: i.e. there is a smooth map
$\varphi_x \colon\bB(0,r_0)\rightarrow \bB(x,r_0)$ with $$\frac 12 \mathrm{eucl}\le \varphi_x ^*g_0\le 2\mathrm{eucl},$$
where $\mathrm{eucl}$ is the Euclidean metric and $\bB(0,r_0)$ the Euclidean ball of radius $r_0$.
Hence there are positive constants $\uptheta$ and $\uplambda$ that only depends on $n$ such that if $\bB$ is a $g_0$ geodesic ball of radius $r(\bB)\le r_0/2$ then 
$$\mu_0(2\bB)\le \uptheta \mu_0(\bB)$$ and
$$\forall \phi\in \cC^1(\bB)\colon
\int_{\bB}\left( \phi-\fint_{\bB} \phi\right)^2d\mu_0\le \uplambda r(\bB)^2\int_{\bB} |d\phi|^2d\mu_0.$$

We can assume that $R_0\le r_0$, hence the above estimation (\ref{BMO}) and the above Poincar\'e inequality implies a BMO (bounded mean oscillation) estimate on the function $\ln (u)$. Moreover, we also get from
(\ref{BMO}) that for all balls $\bB$ with radius $R\leq R_0$, 
$$\int_\bB \abs{d\ln(u)}_{g_0}d\mu_0 \leq C R^{n/2} \left(\int_\bB \abs{d\ln(u)}_{g_0}^2d\mu_0\right)^{1/2} \leq C R^{n-1}.$$
Therefore by John-Nirenberg inequality (cf \cite[Theorem 7.21]{GT}), there exists $\epsilon>0$ such that for any geodesic ball $\bB$ of radius less that $R_0/2$, 
$$\fint_\bB e^{\epsilon \abs{\ln u - (\ln u)_\bB}}d\mu_0 \leq C, $$
where $\ds (\ln u)_\bB = \fint_{\bB} \ln(u) \ d\mu_0$. This implies that
$$\fint_{\bB} u^{\pm \epsilon}\mathrm{d}\mu_0 \le C e^{\pm \epsilon (\ln u)_\bB} .$$
Therefore
$$\left(\fint_{\bB} u^{\epsilon}\mathrm{d}\mu_0\right)\left(\fint_{\bB} u^{-\epsilon}\mathrm{d}\mu_0\right)\le C e^{\epsilon (\ln u)_\bB} e^{- \epsilon (\ln u)_\bB} \leq C .$$

\end{proof}

So far we have only used an estimate of the negative part of the scalar curvature of $g_f$ in $L^{n/2}$. However, we will need the $L^{n/2}$ bound on the positive part of the scalar curvature in order to get the following estimate.

\begin{lem}\label{lem:AInfty2}
Assume the same hypothesis of Theorem \ref{theo:Conformal-AInfty}. 
For all $\epsilon>0$, there exists $C = C(n, g_0, R_0, \delta, \Lambda, \epsilon)>0$ such that for all $g_0$-geodesic balls of radius less than $R_0$,
$$\left(\fint_{\frac 12\bB} u^{\frac{2n}{n-2}}\mathrm{d}\mu_0\right)\le C \left(\fint_{\bB} u^{\epsilon}\mathrm{d}\mu_0\right)^{\frac{2n}{\epsilon(n-2)}}.$$
\end{lem}

\begin{proof}

The proof of this lemma is done in two stages, the first one being a localization of what we did in the proof of Proposition \ref{Integimprove}.
Let $R,r\ge 0$ such that $r+R\le R_0$ and let $\xi$ be a cut-off function with support in $\bB(x,r+R)$ such that 
\begin{itemize}
\item $0\le \xi\le 1$,
\item $|d\xi|_{g_0}\le \frac1r$
\item $\xi=1$ on $\bB(x,R)$
\end{itemize}

We have the equality
\begin{multline*}
\frac{4(n-1)}{n-2}\int_{\bB(x,R+r)}|d(\xi u)|_{g_0}^{2}\mathrm{d}\mu_0\\ 
= \frac{4(n-1)}{n-2}\int_{\bB(x,R+r)}(\Delta_{g_0} u)\,u \xi^{2}\mathrm{d}\mu_0
+ \frac{4(n-1)}{n-2} \int_{\bB(x,R+r)}|d\xi|_{g_0}^{2}u^{2}\mathrm{d}\mu_0\\
=\int_{\bB(x,R+r)}Vu^2 \xi^{2}\mathrm{d}\mu_0 - \int_{\bB(x,R+r)}\scal_{g_0}u^2 \xi^{2}\mathrm{d}\mu_0\\ + \frac{4(n-1)}{n-2} \int_{\bB(x,R+r)}|d\xi|_{g_0}^{2}u^{2}\mathrm{d}\mu_0.
\end{multline*}

We use now the Hebey-Vaugon Sobolev inequality (\ref{eq:HabeyVaugon}) and we get that
\begin{equation}\begin{split}\upalpha(n,2)\left(\int_{\bB(x,R+r)} (\xi u)^{\frac{2n}{n-2}}\mathrm{d}\mu_0\right)^{1-\frac2n}&\le \frac{4(n-1)}{n-2}\int_{\bB(x,R+r)}|d(\xi u)|^{2}\mathrm{d}\mu_0\\
&+B\int_{\bB(x,R+r)}(\xi u)^{2}\mathrm{d}\mu_0
\end{split}
\end{equation}
Hence with H\"older inequality, we get
\begin{equation}\begin{split}
\upalpha(n,2)\left(\int_{\bB(x,R+r)} (\xi u)^{\frac{2n}{n-2}}\mathrm{d}\mu_0\right)^{1-\frac2n} \le \left(B+\lambda_0+\frac{4(n-1)}{n-2} \frac{1}{r^2}\right)\int_{\bB(x,R+r)} u^{2}\mathrm{d}\mu_0\\
+\left(\int_{\bB(R+r)} \left(\scal_{g_f}\right)_+^{\frac n2} \dv_{f}\right)^{\frac 2n}\left(\int_{\bB(x,R+r)} (\xi u)^{\frac{2n}{n-2}}\mathrm{d}\mu_0\right)^{1-\frac2n}.
\end{split}
\end{equation}

We obtain hence that for a constant $C$ that depends only on $g_0,R_0$ and $n$:
\begin{equation}\label{eq:Li-Schoen}
\left(\int_{\bB(x,R)}  u^{\frac{2n}{n-2}}\mathrm{d}\mu_0\right)^{1-\frac2n}\le \frac{C}{\delta r^2}\int_{\bB(x,R+r)} u^2\mathrm{d}\mu_0.
\end{equation}

Therefore,  if $\epsilon\ge 2$, using again H\"older inequality and setting $r=R\leq R_{0}/2$ we get  
$$\left(\int_{\bB}u^{\frac{2n}{n-2}}{\rm d}\mu_0\right)^{\frac{n-2}{n}}\leq \frac{C}{\delta R^2}\left(\int_{2\bB} u^\epsilon{\rm d}\mu_0 \right)^{\frac 2 \epsilon} \mu_0(\bB)^{\frac{\epsilon - 2}{\epsilon}}$$
which gives Lemma \ref{lem:AInfty2}.

 Assume now that $\epsilon<2$ and let us go to the second stage of the proof. We will use a trick of P. Li and R. Schoen \cite[Theorem 2.1]{LiSchoen} in order to obtain the wanted conclusion. Assuming now that $r\le R\le R_0/2$, and using H\"older inequality with exponent
$$\alpha=1+\frac{4}{(2-\epsilon)(n-2)}$$ 
we obtain
$$\left(\int_{\bB(x,R)}  u^{\frac{2n}{n-2}}\mathrm{d}\mu_0\right)^{1-\frac2n}\le \frac{C}{\delta r^2}\left(\int_{B(x,R+r)} u^{\epsilon}\mathrm{d}\mu_0\right)^{1-\frac{1}{\alpha}}\left(\int_{\bB(x,R+r)}  u^{\frac{2n}{n-2}}\mathrm{d}\mu_0\right)^{\frac 1\alpha}.$$

Let us define now $\lambda=\frac{1}{\alpha}\frac{n}{n-2}$,
$$M(R)=\left(\int_{\bB(x,R)}  u^{\frac{2n}{n-2}}\mathrm{d}\mu_0\right)^{1-\frac2n} \quad \mbox{and}\quad \mathbf{I}=\left(\int_{B(x,2R)} u^{\epsilon}\mathrm{d}\mu_0\right)^{1-\frac{1}{\alpha}}.$$
Notice that $\lambda\in (0,1)$ since $\ds \frac 1 \lambda = 1 + \frac{2\epsilon}{n(2-\epsilon)}$. Equation (\ref{eq:Li-Schoen}) shows hence that
$$M(R)\le \frac{C\mathbf{I}}{\delta r^2}\, M(R+r)^\lambda.$$
Iterating this equation with $r_1+r_2+\dots+r_k\le R/2$ leads to
$$M(R)\le\left( \frac{C\mathbf{I}}{\delta}\right)^{1+\lambda+\dots+\lambda^{k-1} }\,\prod_{\ell=1}^k\frac{1}{r_\ell^{2\lambda^{\ell-1}}} M(R+r_1+r_2+\dots +r_k)^{\lambda^k}.$$
If we choose $r_\ell=R2^{-\ell-2}$, with $\beta=\frac{1}{1-\lambda}$, a little bit of arithmetic gives that 
$$M(R)\le C\frac{1}{R^{2\beta}} \mathbf{I}^\beta.$$
It follows that 
$$\left(\fint_{\bB(x,R)} u^{\frac{2n}{n-2}}\mathrm{d}\mu_0\right)\le C \left(\fint_{\bB(x,2R)} u^{\epsilon}\mathrm{d}\mu_0\right)^{\frac{2n}{\epsilon(n-2)}}$$
which concludes the proof of Lemma \ref{lem:AInfty2}.
\end{proof}

Let us now conclude the proof of Theorem \ref{theo:Conformal-AInfty}. Let $\epsilon>0$ be given by Lemma \ref{lem:AInfty1}, and set $p = \frac{2n}{\epsilon(n-2)}$. Let us write $w = u^{\frac{2n}{n-2}} = e^{nf}$. It follows from Lemmas \ref{lem:AInfty1} and \ref{lem:AInfty2} that there exists $C =  C(n, g_0, R_0, \delta, \Lambda, \epsilon)>0$ such that for all $g_0$-balls $\bB$ with radius less than $R_0/2$, 

$$\left(\fint_{\frac 1 2 \bB} wd\mu_0\right)\left(\fint_{\bB} w^{-\frac{1}{p}}d\mu_0\right)^{p}\le C.$$
Hence for any $g_0$-ball $\bB$ with radius less than $R_0/4$:
$$\left(\fint_{ \bB} wd\mu_0\right)\left(\fint_{\bB} w^{-\frac{1}{p}}d\mu_0\right)^{p}\le \uptheta^p \left(\fint_{ \bB} wd\mu_0\right)\left(\fint_{2\bB} w^{-\frac{1}{p}}d\mu_0\right)^{p} \le \uptheta^p C,$$
where $\uptheta$ is the doubling constant for the volume measure $g_0$. This is precisely characterization $ii)$ of $A_\infty$ weights, cf Definition \ref{defi:AInfty}.
\endproof


\section{Strong $A_\infty$ weights and bi-H\"older compactness}\label{sec:SAInfty}

We now introduce so-called \emph{strong $A_\infty$ weights}, also known as \emph{metric doubling measures}, which are special cases of $A_\infty$ weights with strong geometric controls, such as Sobolev, Poincar\'e and isoperimetric inequalities. 

We could not produce examples of sequences of conformal deformations satisfying the hypotheses of Theorem \ref{theo:AInfty} which would not be uniformly strong $A_\infty$. Our Theorem \ref{theo:SAInfty} shows that any conformal deformation is actually uniformly strong $A_\infty$ provided the $L^{n/2}$ norm of the scalar curvature on small balls is pinched enough, which explains this fact.


\subsection{Definitions and basic facts}\label{ssec:SAInfty}

Let $(M,g_0)$ be a closed Riemannian manifold, $\mu_0$ its Riemannian volume. If $g_f=e^{2f}g_0$ is a metric conformal to $g_0$, we continue to denote by $\mu_f$ and by $d_f$ the corresponding Riemannian volume and distance respectively.

\subsubsection{Definitions}

\begin{defi}\label{def:sAInfW}
Let $g_f=e^{2f}g_0$ be a conformal deformation of $g_0$. We say that the volume density $e^{nf}$ is a \emph{strong-$A_\infty$ weight} with respect to the Riemannian metric $g_0$ if there are positive constants $\eta, \uptheta$ such that 
\begin{enumerate}
\item $\ds x\in M, r\le \eta\Rightarrow \mu_f(\bB(x,2r))\le \uptheta\mu_f(\bB(x,r)).$
\item $\ds x,y\in M, d_0(x,y)\le \eta\Rightarrow \frac{d_f(x,y)^n}{\uptheta^{n}} \le \mu_f(\bB(x,d_0(x,y)))\le \uptheta^n d_f(x,y)^n$.
\end{enumerate}
\end{defi}

In above definition, $\bB(x,r)$ is a geodesic ball of radius $r$ for the \emph{fixed} $d_0 = d_{g_0}$ distance. Notice that the strong $A_\infty$ condition is invariant under scaling: if $d\mu_f=e^{nf}d\mu_0$ is a strong-$A_\infty$ weight with constants $(\uptheta,\eta)$ then for every $h>0$, $h^nd\mu_f$ is also a strong-$A_\infty$ weight with same constants $(\uptheta,\eta)$.

These strong $A_\infty$ weights were introduced in \cite{DS}, in relation with the quasiconformal jacobian problem. The jacobian of a quasiconformal map is a strong $A_\infty$ weight, as is any $A_1$-weight (see Definition \ref{defi:AInfty} ii)). A strong $A_\infty$ weight with constants $(\uptheta, \eta)$ is an $A_\infty$-weight  with constants only depending on $\uptheta, \eta$, see Proposition 3.4 of \cite{Semmes1}. Nevertheless, not all $A_\infty$ weights are strong $A_\infty$. Various characterizations of strong $A_\infty$ weights can be found in \cite{DS} and \cite{Semmes1}, together with examples and details about their relationships with the quasiconformal jacobian problem and the bi-Lipschitz embedding problem.

It is easy to see that if the conditions in definition \ref{def:sAInfW} are satisfied for some $\uptheta, \eta$ then for any other $\eta'\ge \eta$ there is a constant $\uptheta'$ depending on  $(\uptheta, \eta)$ and on the geometry of $g_0$ such that the same conditions are satisfied with constant $(\uptheta', \eta')$. We also remark that the doubling condition on $\mu_f$ implies that if $m\in M$ is a point such that $$\rho=d_0(m,x)=d_0(m,y)=\frac 12d_0(x,y)\le \frac\eta 2$$
then the balls $\bB(m,\rho),\ \bB(x,d_0(x,y)),\ \bB(y,d_0(x,y))$ have comparable $\mu_f$ mesure. 

Note that when $d_0(x,y)\le \frac 12 \mathrm{inj}(M,g_0)$, there is a unique point $m\in M$ such that $$\rho=d_0(m,x)=d_0(m,y)=\frac 12d_0(x,y).$$ The ball $\bB(m,\rho)$ will be denoted $\bB_{x,y}$. We obtain thus the following alternative and equivalent definition of strong $A_\infty$ weight.

\begin{defi} If $\uptheta,\eta$  are positive number such that $\eta \le \frac 14 \mathrm{inj}(M,g_0)$, we say that the volume density $e^{nf}$ associated to the metric $g_f=e^{2f}g_0$  is a strong-$A_\infty$ weight (with respect to $g_0)$ with constant $(\uptheta,\eta)$ if 
\begin{enumerate}
\item $\ds x\in M, r\le \eta\Rightarrow \mu_f(\bB(x,2r))\le \uptheta\mu_f(\bB(x,r)).$
\item $\ds x,y\in M, d_0(x,y)\le 2\eta\Rightarrow \frac{d_f(x,y)^n}{\uptheta^{n}} \le \mu_f(\bB_{x,y})\le \uptheta^n d_f(x,y)^n$.
\end{enumerate}
\end{defi} 

It follows from Lemma \ref{lowerd} that the left inequality in property (2) in definition \ref{def:sAInfW} is satisfied for all $A_\infty$ weights. The specificity of strong $A_\infty$ weights relies hence in the other inequality.

\medskip

Let $g_f = e^{2f} g_0$ be a conformal deformation of $g_0$ such that its volume density $e^{n f}$ is a strong-$A_\infty$ weight (with respect to $g_0$) with constants $(\uptheta,\eta)$. We are going to show several basic controls on the geometry of $(M, g_f)$.


\subsubsection{ Geometry of $g_f$ balls.}
In this subsection, we show that the strong $A_\infty$ condition implies that the geodesic ball for $g_f$ and $g_0$ are comparable. 
\begin{center}
Let $B_f(x,r)$ denote the $g_f$-geodesic ball of center $x$ and radius $r$.
\end{center}
In this section and in what follows we use the following notation:
$$R(x,r):=\mu_f(\bB(x,r))^{\frac 1n},$$
where, as before, $\bB(x,r)$ is the $g_0$ geodesic ball of center $x$ and radius $r$. 

\begin{lem}
There exists $\kappa = \kappa(\eta, \uptheta)>0$ such that for any $r\le \eta/\kappa:$
$$B_f(x, R(x, r/\kappa))\subset\bB(x,r)\subset B_f(x,R( x,\kappa r)).$$
\end{lem}
\begin{proof}
The doubling condition implies a reverse doubling estimate  (see \cite[Lemma 2.10]{GS0}): there are positive constants $\alpha,\upomega$ depending only on $\eta$ and $\uptheta$ such that 
\begin{equation}\label{WRatio}
x\in M, 0<s<r\le \eta\Rightarrow\upomega^{-n}\left(\frac{s}{r}\right)^{\frac n\alpha}\le \frac{\mu_f(\bB(x,s)))}{\mu_f(\bB(x,r)))}\le \upomega^n \left(\frac{s}{r}\right)^{n\alpha}.
\end{equation}

The doubling condition also implies that the $\mu_f$-measure of
$g_0$-geodesic spheres is zero (\cite{ColdingMinicozzi, Tessera}), hence the function $r\mapsto R(x,r)$ is continuous and strictly increasing, in particular 
we have $$R(x,s)<R(x,r) \Longleftrightarrow s<r.$$
We have for $d_0(x,y)\le\eta$
\begin{equation}\label{eq:Radius Control}
\uptheta^{-1}\le \frac{d_f(x,y)}{R(x,d_0(x,y))}\le\uptheta .
\end{equation}
and
$$x\in M, 0<s<r\le \eta\Rightarrow \upomega^{-1} \left(\frac{s}{r}\right)^{\frac{1}{\alpha}}\le \frac{R(x,s)}{R(x,r)}\le \upomega \left(\frac{s}{r}\right)^{\alpha}.$$
With $\kappa:=(\upomega \uptheta)^{\frac 1\alpha}$, and assuming without loss of generality that $\kappa >1$, we get for any $r\le \eta/\kappa:$
$$\uptheta R(x,r)\le R(x,\kappa r)\ \mathrm{and}\  R(x,r/\kappa)\le \uptheta^{-1}R(x,r).$$
Since by (\ref{eq:Radius Control}), we have 
$$  B_f(x,R(x,r)/\uptheta)\subset \bB(x,r)\subset B_f(x,\uptheta R(x,r))$$
we get
$$B_f(x, R(x, r/\kappa))\subset\bB(x,r)\subset B_f(x,R( x,\kappa r)).$$
\end{proof}
Note that the two radius $R(x, r/\kappa)$ and $R( x,\kappa r)$ are uniformly comparable:
$$ \omega^{-1} \kappa^{\frac{2}{\alpha}}\le \frac{R( x,\kappa r)}{R( x, r/\kappa)}\le \omega \kappa^{2\alpha}.$$

\subsubsection{ Volume doubling for $(M,d_f,\mu_f)$}

As we mentioned above, if a weight is strong-$A_{\infty}$, it is $A_{\infty}$. We proved that $A_{\infty}$ weights are doubling for the balls of the background metric. We show now that the metric measure space $(M,d_f,\mu_f)$ is doubling (for $g_f$-geodesic balls) with constants only depending on $\eta$ and $\uptheta$.

\begin{lem}
There exists $\delta = \delta(\eta, \uptheta)>0$ and $k = k(\eta, \uptheta)>0$ such that for all $R\in (0, \delta)$,
$$\mu_f(\bB_f(x, 2R)) \leq 2^k \mu_f(\bB_f(x,R)),$$
where $\bB_f(x,R)$ is the $d_f$-geodesic ball of radius $R$ and center $x$.
\end{lem}

\begin{proof}

Let $D$ be the diameter of $(M,g_0)$, then according to \cite[Lemma 3.10]{Carron2016}, the volume doubling condition for $(M, d_0, \mu_f)$ implies that for any $x\in M$ and any $r\le \eta:$
\begin{equation}\label{DiamB1}
\mu_f(M)\le \uptheta^{50+50\frac{D}{r}}\mu_f(\bB(x,r)).
\end{equation}

Hence 
$$R(x,s)\le \left(\uptheta^{-50-50\frac{D}{r}}\mu_f(M)\right)^{\frac 1n}\Rightarrow s\le r.$$
Now we have for $r\le \eta/\kappa:$
\begin{align*}
\mu_f(B_f(x,R(x,r)))&\le \mu_f(\bB(x,{\kappa}r))\\
&\le \uptheta^k \mu_f(\bB(x,2^{-k}r{\kappa}))\\
&\le \uptheta^k \mu_f(B_f(x,R(x,2^{-k}r\kappa^2))
\end{align*}
But
$$R(x,2^{-k}r\kappa^2)\le \omega \left(2^{-k}\kappa^2\right)^\alpha R(x,r)$$
Chosing $k$ such that 
$\omega \left(2^{-k}\kappa^2\right)^\alpha\le \frac 12$ we get
$$\mu_f(B_f(x,R(x,r)))\le  \uptheta^k \mu_f(B_f(x,R(x,r)/2)).$$
Hence with 
$$\delta:=\left(\uptheta^{-50-50\kappa\frac{D}{\eta}}\mu_f(M)\right)^{\frac 1n},$$ 
this eventually implies that 
$$R\le 2\delta\Rightarrow \mu_f(B_f(x,2R))\le 2^k\mu_f(B_f(x,R)).$$
\end{proof}

\subsubsection{Estimate of the diameter of $(M,g_f)$}

We are going to show that the diameter of $(M,g_f)$ is comparable to $\mu_f(M)^{\frac 1n}$.

\begin{lem}\label{lem:DiamStrong} There is a constant $C$ depending only on $g_0,\uptheta,\eta$ such that
$$C^{-1}\left(\mu_f(M)\right)^{\frac 1n} \le\diam(M, g_f)\le C \left(\mu_f(M)\right)^{\frac 1n}.$$
\end{lem}

\begin{proof}
The upper bound for $\diam(M, g_f)$ is an immediate consequence of Lemma \ref{lowerd}. Let us show the lower bound.

It is always possible to assume moreover that $\eta\le D$. Then there is always two points $x,y\in M$ with $d_0(x,y)=\eta$ and we have
$$\diam (M,g_f)\ge d_f(x,y)\ge \uptheta^{-1}\left(\mu_f(\bB(x,\eta))\right)^{\frac 1n},$$
therefore by (\ref{DiamB1}), we get
$$\diam (M,g_f)\ge C(\uptheta,\eta,D)\left(\mu_f(M)\right)^{\frac 1n}.$$
\end{proof}

\subsubsection{Bi-H\"older compactness for sequence of Strong $A_\infty$ conformal metric}

Let $(g_k=e^{2f_k}g_0)_{k\in \N}$ be a sequence of conformal deformations of $g_0$ such that the weights $e^{n f_k}$ are strong-$A_\infty$ weight with uniform constant $(\uptheta,\eta)$. As mentioned earlier, it follows from \cite{Semmes1} that all the $g_k$ are uniformly $A_\infty$. If the volumes of $(M, g_k)$ are uniformly bounded, Proposition \ref{prop:AInfty-Compact} implies then that the sequence of metric spaces $(M, d_{f_k})_{k\in \N}$  has a subsequence which converges in Gromov-Hausdorff topology to some limiting space $(M_\infty, d_ \infty)$. Strong $A_\infty$ weights will give the following stronger convergence result. 

\begin{prop} Let  $(g_k=e^{2f_k}g_0)_{k\in \N}$ be a sequence of Riemannian metrics conformal to $g_0$ such that the weights $e^{n f_k}$ are strong-$A_\infty$ weights with uniform constant $(\uptheta,\eta)$.
\begin{enumerate}[i)]
\item If there are positive constants $v, V$ such that  $\forall k\colon v\le \mu_{f_k}(M)\le V$, then there is a subsequence  that converges in the $\cC^\alpha$ and Gromov-Hausdorff topology to $(M, d_\infty)$, where $d_\infty$ is a distance on $M$ which is bi-H\"older equivalent to $d_0$.
\item If $\lim_{k\to\infty} \mu_{f_k}(M)=0$, then the sequence of metric spaces $\left((M,d_{{f_k}})\right)_k$ converges in the Gromov-Hausdorff topology to a point.
\item If there is a positive constant $v$ such that $\forall k\colon v\le \mu_{f_k}(M)$, then for any $p\in M$, there is a subsequence  of the sequence of pointed metric spaces $\left((M,d_{{f_k}},p)\right)_k$ that converges in the  pointed Gromov-Hausdorff topology to a pointed metric space.
\end{enumerate}
\end{prop}

\begin{proof}

We will only show i); the other cases follow from a straightforward scaling argument. As soon as $e^{nf}$ is a strong $A_\infty$ with constants $\eta, \uptheta$, we have seen that 
if $d_0(x,y)\le \eta$ then
$$\theta^{-1} R(x,d_0(x,y))\le d_f(x,y)\le \theta R(x,d_0(x,y)).$$
Moreover by (\ref{WRatio}) we have
$$ \upomega^{-1}\left(\frac{d_0(x,y)}{\eta}\right)^{\frac 1\alpha} \mu_f(\bB(x,\eta))^{\frac 1n}\le R(x,d_0(x,y))\le \upomega \left(\frac{d_0(x,y)}{\eta}\right)^\alpha \mu_f(\bB(x,\eta))^{\frac 1n}.$$
Hence if $v\le \mu_f(M)\le V$, there is a constant $\uplambda$ depending only on $v,V,\uptheta,\upomega,\eta,\alpha$ such that
$$\frac{1}{\uplambda}\left(d_0(x,y)\right)^{\frac 1\alpha}\le d_f(x,y)\le\uplambda\left(d_0(x,y)\right)^{\alpha}.$$
Hence the identity map $$\mathrm{Id}\colon (M,d_f)\rightarrow (M,d_0)$$ is uniformly bi-H\"older continuous 

Therefore a uniform strong $A_\infty$ estimate for the $(e^{n f_k})_{k\in \N}$ together with a uniform control on the volume $0 < v \leq \mu_{f_k}(M)\leq V$ implies that the sequence $(d_{f_k})$ is precompact in the H\"older topology and for any converging subsequence, the limiting map $d_\infty$ is a distance on $M$ which is uniformly bi-H\"older equivalent to $d_0$.

\end{proof}

\subsubsection{Sobolev, Poincar\'e and isoperimetric inequalities} 

Let us now mention some other important geometric facts about strong $A_\infty$ weights shown in \cite{DS}.

\begin{thm}[David-Semmes (1990)]
Let $(M, g_0)$ be a compact smooth Riemannian manifold, and let $\omega = e^{n f}$ be a strong $A_\infty$ weight with respect to $g_0$ with constants $\eta$ and $\uptheta$, and  let $g_f = e^{2f}g_0 = \omega^{\frac 2 n}g_0$. Then 
\begin{enumerate}[i)]
\item There is a constant $\uptheta=\uptheta(g_0, \eta, \uptheta)$ such that for any $x\in M$ and any $g_f$-geodesic $B(x,r)$ of radius less than $\diam_{g_f}(M)$:
$$\uptheta^{-1}r^n\le  \mu_f(B(x,r))\le \uptheta r^n.$$
\item There is a constant $\upgamma=\upgamma(g_0, \eta, \uptheta)>0$ such that for any smooth domain $\Omega\subset M$ with $\mu_f(\Omega)\le \frac 12 \mu_f(M)$ , we have
$$
\upgamma \mu_f(\Omega)^{\frac{n-1}{n}} \leq  \mu_f(\partial\Omega).$$
\item There is a constant $\uplambda=\uplambda(g_0, \eta, \uptheta)$ such that for any  $g_f$-geodesic $B(x,r)$ of radius less than $\diam_{g_f}(M)$:
and any $\varphi\in\cC^1(B(x,r))$ with $\int_{B(x,r)} \varphi d\mu_f=0$ then
$$\int_{B(x,r)} |\varphi| d\mu_f\le \uplambda r\int_{B(x,r)} |d\varphi|_{g_f} d\mu_f$$
\end{enumerate}\end{thm}
Sobolev inequalities for $(M,g_f)$ are consequences of doubling and Poincar\'e inequality, see \cite{DS}.
\subsection{Q-curvature and Strong $A_\infty$ bounds}

Let us finish this section by presenting the results of S. Brendle \cite{Brendle} and Yi Wang \cite{Wang} which have inspired our analysis of conformal metrics with $L^{n/2}$ pinching of the scalar curvature via $A_\infty$ weights. These works are devoted to the study of conformal deformations with $L^1$ bounds for the Branson $Q$-curvature, which we present now.

\medskip

Let $(M, g)$ be a Riemannian manifold of dimension $n=4$. The Branson $Q$-curvature of $(M, g)$ is defined by
$$
Q_g = \frac 1 {12} \left(\Delta_g \scal_g + \frac 1 4 \scal_g^2 - 3 \abs{E_g}_g\right),
$$
where $E_g$ is the traceless part of the Ricci tensor. This $Q$-curvature has a natural conformal invariance: if $g_f = e^{2f} g_0$, then
\begin{equation}\label{eq:QCurv}
2e^{4 f}Q_{g_f} = P_{g_0} f + 2 Q_{g_0},
\end{equation}
where $P_{g_0}$ is the \emph{Paneitz operator} on $(M, g_0)$, which we shall not define here (see for instance \cite{Brendle}). This equation is the analogous for the $Q$-curvature or the Yamabe equation (\ref{Yamabe}) for the scalar curvature. The existence of a Paneitz operator and a $Q$-curvature satisfying the conformal covariance property (\ref{eq:QCurv}) has been generalized to all even dimension by Fefferman and Graham in \cite{FG02}. Note that the leading term of $P_{g_0}$ is always $(\Delta_{g_0})^{\frac n 2}$, and $P_{g_0} = (\Delta_{g_0})^{\frac n 2}$ when $g_0$ is a flat metric.

\medskip
The following result was shown in \cite{Brendle}, Proposition 1.4.

\begin{thm}[Brendle (2003)]
Let $(\bS^n, g_S)$ be the standard $n$-sphere and let $C, \delta>0$ and $R_0\in (0, {\rm diam}(\bS^n))$ be fixed. Let $(g_k = e^{2 f_k}g_S)_{k\in \N}$ be a sequence of conformal deformations of $g_s$ with constant volume such that:
\begin{itemize}
 \item $\displaystyle \int_{\bS^n} Q_{g_k}^2 d{\rm vol}_{g_k} \leq C$;
 \item for all $x\in \bS^n$, 
 $$\displaystyle \int_{\bB(x,R_0)} \abs{Q_{g_k}}d{\rm vol}_{g_k} \leq \frac 1 2 (n-1)! \sigma_n - \delta,$$
 where $\bB(x,R_0)$ is the $g_S$-geodesic ball of radius $x$ and center $x$.
\end{itemize}

Then $(f_k)_{k\in\N}$ is bounded in $\mathcal H^n = W^{2, n}$ with bounds only depending on $C, \delta, n$ and $R_0$.
\end{thm}

Note that on the round sphere, we have
$$
\int_{\bS^n} \abs{Q_{g_S}}d{\rm vol}_{g_S} = (n-1)! \sigma_n:
$$
the critical constant in Brendle theorem is half of the total $Q$-curvature of the sphere. It was later shown by Wang in \cite{Wang} that under an analogous $L^1$ pinching condition for the $Q$-curvature, conformal deformations of the euclidean metric on $\R^n$ are uniformly strong $A_\infty$ weights.

\begin{thm}[Wang (2013)]
Let $C, \delta>0$ be fixed and let $(g_f = e^{2 f}\mathrm{eucl})_{k\in \N}$ be a conformal deformation of the euclidean plane $(\R^n, \mathrm{eucl})$, which is \emph{normal metrics}, and such that:
\begin{itemize}
 \item $\displaystyle \int_M \abs{Q_{g_f}} d\mu_f \leq C$;
 \item $\displaystyle \int_{M} Q_{g_f}^+d\mu_f \leq \frac 1 2 (n-1)! \sigma_n - \delta.$
\end{itemize}

Then the Riemannian volume density $(e^{n f})_{k\in\N}$ is a strong $A_\infty$ weight with respect to the standard Lebesgue measure on $\R^n$, with constants only depending on $n$, $C$ and $\delta$.
\end{thm}

A \emph{normal} metric on $\R^n$ is a Riemannian metric whose behaviour at infinity is not too pathological, we refer to \cite{Wang} for a precise definition. For instance, any smooth metric whose scalar curvature is non-negative at infinity is normal.

The previous theorem is a reformulation of Theorems 3.2 and 4.1 of \cite{Wang}. The parallel between Brendle and Wang's assumptions and the assumptions of our Theorem \ref{theo:SAInfty} is obvious; indeed our approach was inspired by Wang's proof. Note that the constant $\frac 1 2 (n-1)! \sigma_n$ is optimal to get a strong $A_\infty$ weight, as can easily be shown by considering a sequence of metrics built from removing a ball in $\R^n$ and gluing on its boundary a very long cylinder with a round cap. Such sequence (even smoothed) cannot be uniformly strong $A_\infty$, since it has no uniform isoperimetric inequality.

The proof of Wang can be localized using the same approach as we do in this paper, which gives following result.

\begin{thm}
Let $(M, g_0)$ be a compact smooth manifold, and $C, \delta>0$ and $R_0\in (0, {\rm diam}(M,g_0))$ be fixed. Let $(g_f = e^{2 f}g_0)_{k\in \N}$ be a sequence of conformal deformations of $g_0$, such that:
\begin{itemize}
 \item $\displaystyle \int_M \abs{Q_{g_f}} d\mu_f \leq C$;
 \item for all $x\in M$, 
 $$\int_{\bB(x,R_0)} Q_{g_f}^+d\mu_f \leq \frac 1 2 (n-1)! \sigma_n - \delta.$$
\end{itemize}

Then the Riemannian volume density $(e^{n f})_{k\in\N}$ is a strong $A_\infty$ weight with respect to $g_0$, with constants only depending on $(M,g_0)$, $C$ and $\delta$.
\end{thm}

We will not detail here the proof of this result, which is a combination of our arguments with those of \cite{Wang}.

\section{$L^{n/2}$ pinching for the scalar curvature and strong $A_\infty$ weights}\label{sec:Ln2Strong}

In Section \ref{sec:SAInfty} we explained that a uniform strong $A_\infty$ control coupled with uniform vo\-lu\-me bounds has  important geometric consequences, which were gathered in Corollary \ref{coro:biHolder}.
This section is devoted to the proof  of Theorem \ref{theo:SAInfty}, which is our second main result. In order to motivate the sequel of our study, we start with the proof of Theorem \ref{theo:SAInfty} assuming Theorems \ref{theo:Schrodn2} and \ref{theo:BHS} which will be stated and proved in Sections \ref{ssec:GsSchrOp} and \ref{ssec:W2n2} respectively.

\subsection{Proof of Theorem \ref{theo:SAInfty}}

Let $g_f = e^{2f}g_0$ be a conformal deformation such that for some $\Lambda_0$ small enough (only depending on $g_0$) and for some $R_0>0$, we have
$$
\forall x\in M, \quad \int_{\bB(x,R_0)} \left|\scal_{g_f}\right|^{\frac n2} \dv_{f} \le  \Lambda_0.
$$

By Yamabe's equation (\ref{Yamabe}), writing $e^{2f} = u^{\frac{4}{n-2}}$, we have
$$
\frac{4(n-1)}{n-2} \Delta_{g_0} u = \scal_{g_f} u^{\frac{4}{n-2}}\,u - \scal_{g_0} u.
$$

This is a Schr\"odinger type equation of the form $\Delta_{g_0} u + q_f\ u = 0$, with
$$q_f = \frac{n-2}{4(n-1)}\left(\scal_{g_0} - \scal_{g_f} e^{2f}\right).$$
Moreover, for all $x\in M$,
\begin{eqnarray*}
 \left\|q_f\right\|_{L^{\frac n2}(\bB(x, R_0),g_0)}  &\leq  C(g_0) R_0^2 +\left( \int_{\bB(x, R_0)} \left|\scal_{g_f}\right|^{\frac n2}  d{\rm vol}_{g_f}\right)^{\frac 2n}\\
 &\leq  C(g_0) R_0^2+\Lambda_0.
\end{eqnarray*}
Hence if one chose $\Lambda_0$ to be smaller that the half of the $\upepsilon$ of the Theorem \ref{theo:Schrodn2} , then for $R_0$ small (so that $C(g_0) R_0^2\le \upepsilon/2$) the Schr\"odinger operator $\Delta_{g_0} + q_f$ satisfies the hypothesis of the Theorem \ref{theo:Schrodn2} and
since $u$ is a \emph{positive} smooth solution of $(\Delta_{g_0} + q_f)u = 0$, the lowest eigenvalue of $\Delta_{g_0} + q_f$ is  necessary zero, with eigenfunction $u$. By Theorem \ref{theo:Schrodn2}, we can decompose 
 $$f = h + w,$$ 
 where $h\in W^{2, \frac n 2}(M, g_0)$ and $w\in \cC^\alpha(M, d_0)$, with
$$\sup_{x,y\in M} \frac{|w(x)-w(y)|}{d_0(x,y)^\alpha}\le C(g_0,R_0,\Lambda_0)$$
and
\begin{equation}\label{eq:W2n2}
\|dh\|_{L^n}+\|\Delta h\|_{L^{\frac n2}}\le  C(g_0,R_0, \Lambda_0).
\end{equation}

Let $x_0$ be a chosen base point in $M$. Since $w$ is bounded in the $\cC^\alpha$ topology induced by $d_0$, there exists $C>0$ such that for all $x\in M$,
$$
|w(x) - w(x_0)| \leq C {\rm diam}(M, g_0)^\alpha.
$$

Now, combining (\ref{eq:W2n2}) with Theorem \ref{theo:BHS}, there exists $\eta = \eta(g_0,R_0, \Lambda_0)$ and $\uptheta = \uptheta(g_0,R_0, \Lambda_0)$ such that $e^{nh}$ is a strong $A_\infty$ weight with respect to $(M, g_0)$ with cons\-tants $\eta$ and $\theta$. Since $w - w(x_0)$ is uniformly bounded from above and below, and since being Strong $A_\infty$ weight is a scale invariant condition, this implies that $e^{nf} = e^{n(h + w - w(x_0))} e^{n w(x_0)}$ is also a strong $A_\infty$ weight with constants $\eta$ and $\uptheta$. Moreover, by Proposition \ref{prop:biHolder}, for all $\alpha\in (0,1)$, the distances $d_f$ and $d_0$ are $\alpha$-bi-H\"older with constants only depending on $g_0$, $R_0$, $\Lambda_0$ and $\alpha$, which concludes the proof of Theorem \ref{theo:SAInfty}.

\endproof

\subsection{On the ground state of Schr\"odinger operators with critical potential}
\label{ssec:GsSchrOp}
The purpose of this section is to give a regularity result for the first eigenfunction of a Schr\"odinger operator whose potential has small $L^{n/2}$ norm on small balls, which is the first key step in the proof of Theorem \ref{theo:SAInfty}. The results presented in this paragraph are valid for any compact Riemannian manifold $(M,g)$ and of independent interest.

\begin{thm}\label{theo:Schrodn2} 
Let $(M^n,g)$ be a closed Riemannian manifold of dimension $n>2$. There are  constants $\upepsilon, \alpha>0$ depending only of $(M^n,g)$ such that if $V\in L^{\frac n2}$ satisfies:
\begin{enumerate}
\item the lowest eigenvalue of the Schr\"odinger operator $\Delta-V$ is zero;
\item  there is some  $\rho>0$ such that for any $x\in M$:
$$\int_{\bB(x,\rho)} |V|^{\frac n2}<\upepsilon.$$
\end{enumerate}
Then for all positive $\mathcal C^2$ solution $\varphi>0$ of $\Delta \varphi-V\varphi=0$, there are $f\in W^{2,\frac{ n}{2}}$ and $w\in \cC^\alpha(M)$ such that 
$$\varphi=e^{f+w}.$$
Moreover, writing $\ds \mathbf{I}:=\sup_x \left(\int_{\bB(x,\rho)} |V|^{\frac n2}\right)^{\frac 2n}$ we have
$$\|w\|_{\dot\cC^\alpha}:=\sup_{x,y\in M} \frac{|w(x)-w(y)|}{d_0(x,y)^\alpha}\le C(g,\rho,\mathbf{I})$$
and
$$\|df\|_{L^n}+\|\Delta f\|_{L^{\frac n2}}\le  C(g,\rho)\mathbf{I}.$$
\end{thm}

This theorem is a key step in showing that conformal deformations whose scalar curvature has small $L^{n/2}$ norm on small balls are strong $A_\infty$ weights.

\medskip

In order to prove Theorem \ref{theo:Schrodn2}, we need some preliminary tools. Let us start by recalling the following classical consequence of the fixed point theorem. 

\begin{prop}\label{Solu} Let $(X,\|\bullet\|)$ be a Banach space, $\kappa\in (0,1)$. If $S\colon \cB(0,\rho)\rightarrow  X$ is a $\kappa$-Lipschitz map with
$$\|S(0)\|+\kappa\rho<\rho$$ then there is a unique $x\in  \cB(0,\rho)$ such that $S(x)=x$. Moreover we have $\ds \|x\|\le \frac{\|S(0)\|}{1-\kappa}$.
\end{prop}

Let $A(g)$ be the best constant in the following functional inequality on the closed Riemannian manifold $(M,g)$:
\begin{equation}\label{dDelta}
\forall \varphi\in \cC^\infty(M)\colon \|d\varphi\|_{L^n}\le A(g)\|\Delta \varphi\|_{L^{\frac n 2}}.
\end{equation}

\begin{prop}
Let $(M,g)$ be a compact smooth Riemannian manifold such that 
${\rm Ric}_{g}\ge -(n-1)k^2 g$, then the best constant $A(g)$ in (\ref{dDelta}) satisfies
$$A(g)\le e^{c(n)(1+k \diam (M,g))}\frac{\diam (M,g)}{\mathrm{vol}^{\frac 1n}(M,g)}.$$
\end{prop}

The above estimate on $A(g)$ can be shown using a lower bound on the Ricci curvature, and upper bounds  on the diameter and volume of $(M,g)$; it follows from results of S. Gallot  \cite{Gallot} and of D. Bakry \cite[Th\'eor\`eme 4.1, Lemme 4.2 and Corollaire 4.3]{Bakry}. We will not show it here, since we will not need this expression.

Let us now show the following property of the ground state of the Schr\"odinger operator $\Delta-V$.

\begin{prop}\label{Gstate}Let $V\in L^{\infty}(M)\cap L^{\frac n 2}(M)$ and let $\lambda_0$ be the bottom of the spectrum of the Schr\"odinger operator $\Delta-V$. If  
$$8A(g)^2\|V\|_{L^{\frac n2}}<1$$ then for  any $\varphi>0$ satisfying $\Delta\varphi-V\varphi=\lambda_0 \varphi$, we have
$$\|d\log \varphi\|_{L^n}\le 2A(g)\|V\|_{L^{\frac n2}}.$$
\end{prop}

\proof We remark that if $\varphi=e^u$ satisfies $$\Delta\varphi-V\varphi=\lambda_0 \varphi,$$ then we have
$$\Delta u-|du|^2=V+\lambda_0.$$
Moreover up to a scaling constant, there is a unique non trivial solution to the equation $\Delta\varphi-V\varphi=\lambda_0 \varphi$.
We introduce the Banach space
$$W^{1,n}_\perp:=\{v\in W^{1,n}(M), \int_M v=0\}$$ endowed with the norm
$$\|v\|_\perp=\|dv\|_{L^n}.$$
We introduce the operator $S\colon W^{1,n}_\perp\rightarrow W^{1,n}_\perp$ defined by
$$S(v)=\Delta^{-1} V+\Delta^{-1}|dv|^2.$$
where for all $f\in L^{\frac n2}$, we have denote 
$u:=\Delta ^{-1} f$ the solution of the equation
$$\begin{cases} \Delta u=f-f_M&\\
\int_M u=0.&
\end{cases}$$
By definition of $A(g)$, we have
$$\|S(0)\|_\perp\le A(g)\|V\|_{L^{\frac n2}},$$ and
$$\|S(u)-S(w)\|_\perp\le A(g)\|u-w\|_\perp \left(\|u\|_\perp+\|w\|_\perp\right).$$
As a consequence the restriction of $S$ to the ball $\cB(0,\rho)\subset W^{1,n}_\perp$, is $2A(g)\rho$-Lipschitz and if
\begin{equation}\label{ineq}
A(g)\|V\|_{L^{\frac n2}}+2A(g)\rho^2<\rho
\end{equation} then there is a unique $v\in \bB(0,\rho)$ with $S(v)=v$ and 
$$\|v\|_\perp\le A(g)\frac{\|V\|_{L^{\frac n2}}}{1-2A(g)\rho}.$$
Choosing $\rho=(4A(g))^{-1}$, the condition (\ref{ineq}) is satisfied as soon as
$$A(g)\|V\|_{L^{\frac n2}}< (8A(g))^{-1}.$$ 
We have obtained hence a function $v\in W^{1,n}_\perp$ such that 
$S(v)=v$ and 
$\|v\|_\perp\le 2A(g)\|V\|_{L^{\frac n2}}.$
Therefore there is a constant $c$ such that 
$$\Delta v-|dv|^2=V+c,$$
therefore
$$\Delta e^v-Ve^v=ce^v.$$
Since the only positive eigenfunction of $\Delta - V$ is associated to $\lambda_0$, we have $c=\lambda_0$, $e^v = k\varphi$ for some $k>0$ and $\ds \norm{d\log \varphi}_{L^n} = \norm{v}_\perp \leq 2A(g)\norm{V}_{L^{n/2}}$.
\endproof

\medskip

Let us now define the Sobolev constant $\beta = \beta(g_0)>0$ of $(M,g_0)$ to be the best constant in the inequality
\begin{equation}\label{Sobolev}
\forall \varphi\in \cC^\infty(M)\colon \ \beta\|\varphi\|_{L^{\frac{2n}{n-2}}}^2\le \|d \varphi\|^2_{L^{2}}+ \| \varphi\|^2_{L^{2}}.
\end{equation}

Note that this Sobolev inequality implies that 
$$\beta\le \vol(M,g_0)^{\frac 2n}.$$
We also have that for any $\delta\in(0,1)$:
$$\forall \varphi\in \cC^\infty(M)\colon \ \delta^2\beta\|\varphi\|_{L^{\frac{2n}{n-2}}}^2\le \int_M\left[ |d \varphi|^2_{g_0}+ \delta^2 |\varphi|^2\right]d\mu_0.$$
Hence we can easily prove the following lemma:
\begin{lem}\label{Ln2 positive} If $W\in L^{\frac n2}$ satisfies 
$\|W\|_{L^{\frac n2}}\le \beta\delta^2$,  then the Schr\"odinger operator $\Delta+\delta^2+W$ is non negative.
\end{lem}
\begin{proof} The H\"older inequality  implies that if $\varphi\in \cC^1(M)$:
$$
-\int_M W\varphi^2\le \|W\|_{L^{\frac n2}}\|\varphi\|_{L^{\frac{2n}{n-2}}}^2\le \delta^2\beta\|\varphi\|_{L^{\frac{2n}{n-2}}}^2\le \int_M\left[ |d \varphi|^2_{g_0}+ \delta^2 |\varphi|^2\right]d\mu_0
$$
\end{proof}

\begin{lem}\label{GS} Assume that $x_0\in M$ and $r_0$ are chosen so that 
$$\vol \bB(x_0,r_0)\le \left(\frac{\beta}{2}\right)^{\frac n2}.$$
If $q\in L^{\frac n 2}$ is such that $$\mathrm{supp} (q)\subset \bB(x_0,r_0)\ \mathrm{and}\   \|q\|_{L^{\frac n2}}\le \beta$$ then there is a constant $c_0$ satisfying :
$$|c_0|\le \frac{2}{\beta} \|q\|_{L^{\frac n2}}$$ such that the lowest eigenvalue of the  Schr\"odinger operator
$$\Delta+q+c_0\un_{M\setminus \bB(x_0,r_0)}$$ is zero.
\end{lem}

\begin{proof} Let $\lambda_0(c)$ be the lowest eigenvalue of the Schr\"odinger operator 
$$\Delta+q+c\un_{M\setminus \bB(x_0,r_0)}.$$
Using the constant function, we get that  $\lambda_0(c)\le 0$ as soon as
$$\int_M q+c\vol_{g_0}(M\setminus \bB(x_0,r_0))\le 0.$$
Moreover, $\vol \bB(x_0,r_0)\le\frac{ \vol (M,g_0)}{2},$ hence for
$c_-=-\frac{2}{\vol(M,g_0)}\int_M |q|$, we have
$\lambda_0(c_-)\le 0$.

Using Lemma \ref{Ln2 positive} for $W = q - c\un_{\bB(x_0,r_0)}$, we also see that if $c>0$ satisfies 
$$
\|q-c\|_{L^{\frac n2}(\bB(x_0,r_0))}\le \beta c
$$ then 
$\lambda_0(c)\ge 0$.
For instance for $c_+=\frac{2}{\beta}\|q\|_{L^{\frac n2}}$, we have
$\lambda_0(c_+)\ge 0$. Hence there is a $c_0\in [c_-,c_+]$ such that $\lambda_0(c_0)=0$. 

Moreover, by H\"older inequality we have
$$|c_-|\le \frac{2}{\vol_{g_0} (M)}(\vol_{g_0} (M))^{1-\frac 2 n}\|q\|_{L^{\frac n2}}\le\frac{2}{\vol^{\frac{2}{n}}_{g_0} (M)}\|q\|_{L^{\frac n2}}\le \frac{2}{\beta} \|q\|_{L^{\frac n2}}.$$
Therefore $\displaystyle |c_0|\leq \frac{2}{\beta} \|q\|_{L^{\frac n2}}$, which concludes the proof of Lemma \ref{GS}.

\end{proof}

We can now prove Theorem \ref{theo:Schrodn2}.

\begin{proof}[Proof of Theorem \ref{theo:Schrodn2}]
We will use the constants $A(g_0)$ and $\beta$ given by the estimate (\ref{dDelta}) and the Sobolev inequality (\ref{Sobolev}). We can always assume that $\rho$ also satisfies 
$$\forall x\in M\colon\, \vol_{g_0}\bB(x,\rho)\le \left(\frac{\beta}{2}\right)^{\frac n 2}.$$
We choose a finite cover 
$$M=\cup_{i=1}^N \bB(x_i,\rho/2)$$ such that for any $i\not=j:$ $\bB(x_i,\rho/4)\cap \bB(x_j,\rho/4)=\emptyset$.
If $$\mathbf{I}:=\sup_x \left(\int_{\bB(x,\rho)} |V|^{\frac n2}\right)^{\frac 2n}\le \frac{\beta}{2}$$
then for any $i\in\{1,\dots,N\}$, using Lemma \ref{GS}, we find $c_i$ such that with $$q_i=V\un_{\bB(x_i,\rho)}+c_i\un_{M\setminus \bB(x_i,\rho)},$$  the lowest eigenvalue of the  Schr\"odinger operator $\Delta-q_i$ is zero.
We have
$$\|q_i\|_{L^{\frac n 2}}\le \frac{3 \vol_{g_0}(M)^{\frac 2n}}{\beta}\,\|V\|_{L^{\frac n 2}(\bB(x_i,\rho))}.$$
Hence using Proposition \ref{Gstate}, we know that if
$$24\frac{ \vol_{g_0}(M)^{\frac 2n}}{\beta}\,\mathbf{I}\,A^2(g_0)<1$$ then we find $\varphi_i\in \ker(\Delta-q_i)$, $\varphi>0$, with
$$\| d\log\varphi_i\|_{L^n}\le C(g_0)\,\mathbf{I}.$$
On the other hand,
$$\Delta \log\varphi_i=q_i+| d\log\varphi_i|^2$$
hence
$$\| \Delta\log\varphi_i\|_{L^{\frac n2}}\le C(g_0) \mathbf{I}.$$
We choose $\upepsilon$ such that  $24\frac{ \vol_{g_0}(M)^{\frac 2n}}{\beta}\,A^2(g_0)\,\upepsilon=\frac12$.
Let $\varphi : M\to (0,+\infty)$ be a global solution of $\Delta\varphi = V\varphi$. For all $i\in \{1,..., N\}$, we define $\psi_i\colon \bB(x_i,\rho)\rightarrow \R$ by
$$\varphi=\varphi_i\psi_i.$$ Since on the ball $\bB(x_i,\rho)$, we also have $\Delta\varphi_i = V\varphi_i$, we get
$$\Delta \psi_i-2\langle d\log \varphi_i,d\psi_i\rangle=0.$$
Using the work of O. A. Ladyzhenskaya and N. N. Ural'tseva \cite{LU85} and of M.V. Safonov \cite{Safo}, we know that 
$\psi_i$ is H\"older continuous.
There is some $\alpha=\alpha(g_0, \rho,\upepsilon)$ such that 
$$\|\psi_i\|_{\dot\cC^\alpha(\bB(x_i,3\rho/4))}:=\sup_{x,y\in \bB(x_i,3\rho/4))} \frac{|\psi_i(x)-\psi_i(y)|}{d_0(x,y)^\alpha}\le C(g_0,\rho, \mathbf{I})\inf_{x\in \bB(x_i,3\rho/4))} \psi_i(x).$$
We let $\delta_i=\inf_{\bB(x_i,3\rho/4)} \psi_i$ and $\psi_i=\delta_ie^{w_i}$ and $\ds \varphi_i=\frac{e^{f_i}}{\delta_i}$ and $\overline{f_i}=\fint_{\bB(x_i,3\rho/4)}f_i$.
We have
$$\|w_i\|_{L^\infty(\bB(x_i,3\rho/4))}+\|w_i\|_{\dot\cC^\alpha(\bB(x_i,3\rho/4))}\le C(g_0,\rho, \mathbf{I})$$ and 
$$\| d f_i\|_{L^n}+\| \Delta f_i\|_{L^{\frac n2}}\le C(g_0) \mathbf{I}.$$
Moreover, by Sobolev inequality
$$\|  f_i-\overline{f_i}\|_{L^n}\le C(g_0) \mathbf{I}.$$

Note that on all ball $\bB(x_i, \rho)$, we have
\begin{equation}\label{eq:DecompPhi}
\varphi = \varphi_i \psi_i = e^{w_i} e^{f_i}.
\end{equation}

Let us fix now a partition of unity $(\chi_i)_{i\in \{1,...,N\}}$ with $\mathrm{supp}\chi_i\subset \bB(x_i,3\rho/4)$ and $\chi_i=1$ on $\bB(x_i,\rho/4)$. This partition of unity only depends on $(M, g_0)$.
Let us define
$$f:=\sum_i \chi_i \left(f_i-\overline{f_i}\right) \quad \mbox{and} \quad w = \sum_i \chi_i(w_i + \overline{f_i}).$$
The functions $f$ and $w$ are well defined on $M$ and it follows from (\ref{eq:DecompPhi}) that $\varphi=e^fe^w$.

We have moreover
$\| d\log f\|_{L^n}+\| \Delta f\|_{L^{\frac n2}}\le C(g_0) \mathbf{I}.$
We also have $\log\varphi=f+w$ and 
$$\Delta \log\varphi-|d\log\varphi|^2=V$$
Hence by Stokes formula,
$$\int_M |d\log\varphi|^2d\mu_0\le \int_M |V|d\mu_0\le \left(\mu_0(M)\right)^{1-\frac 2n}\, \|V\|_{L^{\frac n2}}\le C(g_0)\mathbf{I}.$$
Using the Poincar\'e inequality, one conclude that for $c=\fint_M \log(\varphi)d\mu_0$:
$$\|f+w-c\|_{L^2}\le C(g_0)\mathbf{I}.$$
On each ball $\bB(x_i,\rho/4)$, the oscillation of $w_i=w-\overline{f_i}$ are controlled hence one get that 

$$\left|\overline{f_i}-c\right|\le C(g_0)\|f_i-c\|_{L^2(\bB(x_i,\rho/4))}\le C(g_0,\mathbf{I}).$$
Hence $w=\sum \chi_i\left(w_i+\overline{f_i}\right)$ is bounded in $\dot\cC^\alpha(M)$.
\end{proof}
\subsection{Strong $A_\infty$ weights from $W^{1, n}$ bounds}\label{ssec:W2n2}

In this section, we explain how the proof of Theorem 3.1 of \cite{BHS2} can be adapted for closed Riemannian manifolds, which gives the following result. This is the second key step in the proof of Theorem \ref{theo:SAInfty}.

\begin{thm}\label{theo:BHS} Let $(M,g_0)$ be a closed Riemannian manifold. If $f\colon M\rightarrow \R$ is such that
$$\|d f\|_{L^{n}}\le \mathbf{I},$$ 
for some $\mathbf{I}$. Then the weight $e^{nf}$ is a strong $A_\infty$ weight with respect to $g_0$, with constants $(\uptheta,\eta)$
where $\eta$ depends only on $g_0$ and $\uptheta$ depends only on $g_0$ and $\mathbf{I}$.
\end{thm}

Note that, unlike Theorem 3.1 of \cite{BHS2}, our statement and our proof are scale invariant. 

\subsubsection{Presentation of the setting}

The estimates will depend on several geometric and analytic estimates. We continue to denote by $\mu_0$ the Riemannian volume measure and by $d_0$ the Riemannian distance associated to $g_0$. There is some $\eta>0$ and positive constants $\uptheta,\upgamma,C$ (that can be chosen depending only on the dimension $n$)  such that for any ball $\bB\subset M$ with radius $r(\bB)\le 2\eta$ we have  

\begin{enumerate}[i)]
\item as already mentioned in (\ref{eq:rel-vol-rad-mu0}), $\uptheta^{-1 }r(\bB)^n\le \mu_0(\bB)\le\uptheta r(\bB)^n$;
\item as shown in \cite{GT}, Lemma 7.16, for any Lipschitz function $\varphi\colon \bB\rightarrow \R$ and any $x\in \bB$:
$$
\left|\varphi(x)-\varphi_\bB\right|\le C\int_{\bB} \frac{ |d\varphi|_{g_0}(z)}{d_0^{n-1}(x,z)}d\mu_0(z);
$$
\item the following Adams-Moser-Trudinger inequality: for any Lipschitz function $\varphi\colon \bB\rightarrow \R$ and any $x\in \bB$:
$$\fint_{\bB} e^{\left(\frac{\gamma \left|\varphi(x)-\varphi_\bB\right|}{\|d\varphi\|_{L^{n}(\bB)}}\right)^{\frac{n}{n-1}}}d\mu_0(x)\le C .$$

For the Euclidean space $\R^n$, this inequality is due to N. Trudinger \cite{Trudinger67}. According to \cite[Lemma 7.13]{GT}, the first two geometric control i-ii) implies the Trudinger inequality.  
\end{enumerate}

\item We will use several times the following integration by parts formula
\begin{eqnarray*}
\int_{\bB(x,R)\setminus \bB(x,r)}\varphi(d_0(x,y))f(y)d\mu_0(y)=\varphi(R)\int_{\bB(x,R)}f(y)d\mu_0(y)\\
-\varphi(r)\int_{\bB(x,r)}f(y)d\mu_0(y)-\int_r^R\varphi'(r)\int_{\bB(x,r)}f(y)d\mu_0(y)dr.
\end{eqnarray*}

\subsubsection{Estimate on $\mu_f(\bB)$} 

We are going to prove that $\mu_f$ is an $A_p$ weights for every $p>1$ using characterization ii) in Definition \ref{defi:AInfty}. Let $\bB(o,\rho)\subset M$ be such that its radius satisfies $r(\bB)\le \eta$. We have for any $x\in \bB$ and any $\lambda\in  \R$:
$$\lambda |f(x)-f_\bB|\le  \frac{n-1}{n}\left(\frac{\gamma \left|f(x)-f_\bB\right|}{\|df\|_{L^{n}(\bB)}}\right)^{\frac{n}{n-1}}+\frac 1 n \left(\frac{|\lambda|}{\gamma}\| d f\|_{L^{ n}(\bB)}\right)^{n}.$$
Hence, using Adams-Moser-Trudinger inequality, we get
\begin{align*}
\fint_{\bB}e^{\lambda f}d\mu_0&\le e^{\lambda f_\bB}e^{\left(\frac{|\lambda|}{c}\| d f\|_{L^{n}(\bB)}\right)^{n}}\fint_\bB \exp\left(\left(\frac{\gamma \left|f(x)-f_\bB\right|}{\|df\|_{L^{n}(\bB)}}\right)^{\frac{n}{n-1}}\right)d\mu_0(x)\\
&\le C(g_0, |\lambda|, \mathbf{I}) e^{\lambda f_\bB} 
 \end{align*}
Using the Cauchy-Schwarz inequality, this implies
 $$1\le \fint_{\bB}e^{\lambda f}d\mu_0\times\fint_{\bB}e^{-\lambda f}d\mu_0\le C(g_0, |\lambda|,\mathbf{I}) \fint_{\bB}e^{\pm\lambda f}d\mu_0  e^{\lambda f_\bB},$$
and we eventually get 
\begin{equation}\label{AverLambdaf}
\frac{e^{\lambda f_\bB}}{C(g_0, |\lambda|,\mathbf{I})}\le \fint_{\bB}e^{\lambda f}d\mu_0 \le C(g_0, |\lambda|,\mathbf{I}) e^{\lambda f_\bB}.
\end{equation}

This easily implies that for any $p>1$ there is a constant $C$ depending on $g_0,\mathbf{I},p$ such that 
\begin{equation}\label{eq: W1n implies Ap}
\fint_{\bB }e^{nf}d\mu_0\left(\fint_{\bB }e^{-\frac{n}{p-1}f}d\mu_0\right)^{p-1}\le C.
\end{equation}
Hence, $e^{nf}$ is an $A_p$ weight, in particular  $d\mu_f$ is a doubling measure.
Moreover, by Lemma \ref{lowerd}, this implies that there exists a constant $C = C(g_0,\mathbf{I})>0$ such that 
$$d_0(x,y)\le 10\eta\Rightarrow d_f(x,y)\le C \left(\mu_f(\bB_{x,y})\right)^{\frac 1 n}.$$

\subsubsection{Distance estimate}

 To show that $d\mu_f$ is a strong $A_\infty$ weight, we are only left with showing that for all $x,y\in M$ and for some $C>0$,
$$d_{g_f}(x,y)\geq C\mu_f(\bB_{x,y})^{\frac 1 n}.$$ 
In this purpose, we will need the following lemma:
 \begin{lem}\label{H1} For any $\upepsilon>0$ there is a constant $C$ depending on $g_0, \mathbf{I}$  and $\upepsilon>0$ such that for any ball $\bB\subset M$ with $r(\bB)\le \eta$ there is  a set $E\subset \bB$ such $E\subset \bigcup_\alpha \bB_\alpha$ with
 $$\sum_\alpha r(\bB_\alpha)\le \upepsilon r(\bB)$$ and
 $$x\in \bB\setminus E\Rightarrow |f(x)-f_\bB|\le C.$$
 \end{lem}
 \proof
 Assume that $x\in \bB$ is such that for any $r\in (0,2r(\bB))$ we have
 $$\int_{\bB(x,r)}|df|^{n}\le \lambda\frac{r}{r(\bB)}$$ then using 
 $V(x,r):=\int_{\bB(x,r)}|d f(y)|d\mu_0(y)$ and one gets 
\begin{align*}
|f(x)-f_\bB|&\le \int_{\bB(x,2r(\bB))}\frac{C}{d^{n-1}_0(x,y)}|d f(y)|d\mu_0(y)\\
&\le\frac{C}{(2r(\bB))^{n-1}}V(x,2r(\bB))+\int_0^{2r(\bB)}\frac{(n-1)C}{r^{n}}V(x,r)dr
\end{align*}
With H\"older inequality, one has 
$$V(x,r)\le \uptheta^{1-\frac1n} r^{n-1}\left(\int_{\bB(x,r)}|d f|^{n}\right)^{1/n}\le \uptheta^{1-\frac1n} r^{n-1}\left(\lambda\frac{r}{r(\bB)}\right)^{1/n}.$$
 Hence there is a constant depending only of $g_0$ such that 
 $$|f(x)-f_\bB|\le C \lambda^{\frac 1n}.$$
 If now 
 $E$ is the subset of  point $x$ where there is a $r_x\in (0,2r(\bB))$ such that
 $$\int_{\bB(x,r_x)}|d f|^{n}>\lambda\frac{r_x}{r(\bB)}$$
 Then using Vitali covering lemma, one gets a covering
 $E\subset \cup \bB_\alpha$ where
 $\bB_\alpha=5\bB(x_\alpha,r_\alpha)$ such that 
$$ \int_{\bB(x_\alpha,r_\alpha)}|d f|^{n}>\lambda\frac{r_\alpha}{r(\bB)}$$
and the balls $\bB(x_\alpha,r_\alpha)$ are disjoints. In particular this implies that
$$\sum r(\bB_\alpha)=\sum 5r_\alpha\le \frac{5r(\bB)}{\lambda}\sum \int_{\bB(x_\alpha,r_\alpha)}|d f|^{n}\le \frac{5r(\bB)}{\lambda} \int_M|d f|^{n}.$$
Choosing 
$\lambda=\frac{5}{\upepsilon\mathbf{I}^n}$, we get the desired result.
 \endproof
 
 Let now $x,y\in M$ such that $\rho:=d_0(x,y)\le 10\eta$ and let $\gamma\colon [0,1]\rightarrow M$ be a $\cC^1$ curve joining $x$ to $y$:
 $$\gamma(0)=x\ \mathrm{and}\ \gamma(1)=y$$
 Let $\tau>0$ be the first time with $d_0(x,\gamma(\tau))=\rho$ and let $\bar y=\gamma(\tau)$, $\bB:=\bB_{x,\bar y}=\bB(m,\rho/2)$. We have for any
 $t\in [0,\tau], \gamma(t)\in 3\bB$.
 Let $c(t)=\exp_m(tu)$ be the unit speed $g_0$ geodesic joining $x$ to $\bar y$ and $L:=c([-3\rho/2,3\rho/2])$ it is a diameter of $3\bB$ and 
 $\Pi\colon 3\bB\rightarrow L$ the projection onto $L$ associated to the Euclidean metric $d_{\mathrm{eucl}}$ induced by $\exp_m$.
 We know that $\Pi$ is $1$-Lipschitz for the Euclidean distance and for $x,y\in L$ : $d_{\mathrm{eucl}}(x,y)=d_{0}(x,y)$. Hence $\Pi$ is $\kappa$-Lipschitz for the distance $d_0$.
According to Lemma \ref{H1} with  
\begin{equation}\label{choice}
\upepsilon=\frac{1}{6\kappa}
\end{equation}
we can find $E\subset 3\bB$ such that 
 $E\subset \bigcup_\alpha \bB_\alpha$ with
 $$\sum_\alpha r(\bB_\alpha)\le \upepsilon\frac 32 \rho$$ and
 $$x\in 3\bB\setminus E\Rightarrow |f(x)-f_\bB|\le C.$$

 We have
 \begin{align*}
 L_{g_f}(\gamma)&\ge L_{g_f}(\left.\gamma\right|_{[0,\tau]})\ge L_{g_f}(\left.\gamma\right|_{[0,\tau]}\setminus E)\\
 &\ge e^{-C}e^{f_\bB}L_{g_0}\left(\left.\gamma\right|_{[0,\tau]}\setminus E\right)\\
 &\ge e^{-C}e^{f_\bB}\kappa^{-1}L_{\mathrm{eucl}}\left(\left.\gamma\right|_{[0,\tau]}\setminus E\right)\\
 &\ge e^{-C}e^{f_\bB}\kappa^{-1}L_{\mathrm{eucl}}\left(\Pi\left(\left.\gamma\right|_{[0,\tau]}\setminus E\right)\right))
 \end{align*}
 Since $\Pi\left(\left.\gamma\right|_{[0,\tau]}\right)$ contains the geodesic segment joining $x$ to $\bar y$,
 $$L_{\mathrm{eucl}}\left(\Pi\left(\left.\gamma\right|_{[0,\tau]}\right)\right)\ge \rho.$$
 whereas
 $$L_{\mathrm{eucl}}\left(\Pi(E)\right)\le L_{\mathrm{eucl}}\left(\Pi(\cup \bB_\alpha)\right)\le \kappa\sum_\alpha 2r(\bB_\alpha)\le 2\kappa \upepsilon \frac{3}{2} \rho.$$
 With the choice of $\upepsilon$ given by (\ref{choice}), we get
 $$L_{\mathrm{eucl}}\left(\Pi(E)\right)\le \frac 12\rho$$ hence
 $$L_{\mathrm{eucl}}\left(\Pi\left(\left.\gamma\right|_{[0,\tau]}\setminus E\right)\right))\ge \frac 12\rho.$$ 
 Eventually using (\ref{AverLambdaf}) we obtain:
 $$ L_{g_f}(\gamma)\ge C e^{f_\bB} r\rho\ge C \left(\mu_f(\bB_{x,y})\right)^{\frac 1n}.$$
 
This concludes the proof of Theorem \ref{theo:BHS}.
 
 \medskip

Let us now show that the previous proof implies that the distances $d_f$ and $d_0$ are $\alpha$-bi-H\"older for all $ \alpha\in (0,1)$.

\begin{prop}\label{prop:biHolder}
 Let $(M,g_0)$ be a closed Riemannian manifold. If $f\colon M\rightarrow \R$ is such that
$$\|d f\|_{L^{n}}\le \mathbf{I}$$  
Then for each $\alpha\in (0,1)$ there is a constant $C(g_0, \mathbf{I}, \alpha)$ such that for any $x,y\in M$
$$\frac{\mu_f(M)^{\frac{1}{n}}}{C(g_0, \mathbf{I}, \alpha)} \, d^{\frac{1}{\alpha}}_0(x,y)\le d_f(x,y)\le C(g_0, \mathbf{I}, \alpha)\, \mu_f(M)^{\frac1 n}\,d^\alpha_0(x,y).$$
\end{prop}
\begin{proof}
 We have shown in Theorem \ref{theo:BHS} that the two weights $e^{nf}$ and $e^{-nf}$ are Strong $A_\infty$ weight. In particular, there is a constant $C$ depending on $g_0$ and $\mathbf{I}$ such that for all $x,y\in M$:
\begin{align*}
\frac{1}{C}\,d_f^n(x,y)&\le \mu_f\left(\bB_{x,y}\right)\le C\, d^n_f(x,y)\\
\frac{1}{C}\,d_{-f}^n(x,y)&\le \mu_{-f}\left(\bB_{x,y}\right)\le C\, d^n_{-f}(x,y)
\end{align*}

The estimate (\ref{AverLambdaf}) implies that for each $\uplambda\in \R$ there is a constant $C'$ depending on $g_0$, $\uplambda$ and $\mathbf{I}$ such that for any $x,y$:
$$\frac{1}{C'}\, d^n_0(x,y)\le\left(\int_{\bB_{x,y}} e^{\uplambda f}d\mu_0\right)\left(\int_{\bB_{x,y}} e^{-\uplambda f}d\mu_0\right)\le C'\, d^n_0(x,y).$$
Hence, there is a constant $\Lambda$ such that for any $x,y\in M$:
\begin{equation}\label{dfdmf}
\frac{1}{\Lambda}\,  d^{2n}_0(x,y)\le d_f^n(x,y)\,d_{-f}^n(x,y)\le \Lambda\, d^{2n}_0(x,y).\end{equation}
Notice that this implies that for some constant $C$:
$$\frac{1}{C} \le \mu_{f}(M)\mu_{-f}(M)\le C.$$
But with the H\"older inequality, we get that for any $p>1$ and $x,y\in M$:
$$d_f^n(x,y)\le C d_0^n(x,y)\left(\fint_{\bB_{x,y}} e^{\frac{p}{p-1} f}d\mu_0\right)^{1-\frac{1}{p}},$$ and again the estimate (\ref{AverLambdaf}) implies that
$$\left(\fint_{\bB_{x,y}} e^{\frac{p}{p-1} f}d\mu_0\right)^{1-\frac{1}{p}}\le C\, \frac{\mu_f(\bB_{x,y})}{d_0^n(x,y)}\le  C\, \frac{\mu_f(M)}{d_0^n(x,y)}.$$
So that we get that for any $x,y\in M$ and any $p>1$, 
$$
d_f^n(x,y)\le C\, d^{\frac{n}{p}}_0(x,y) \mu_f(M)
$$
which is the right-hand side inequality of Proposition \ref{prop:biHolder}. Moreover, since
$$
d_{-f}^n(x,y)\le C\, d^{\frac{n}{p}}_0(x,y)\mu_{-f}(M)\le  C\, d^{\frac{n}{p}}_0(x,y)\left(\mu_f(M)\right)^{-1}, 
$$
using the comparison (\ref{dfdmf}) we get:
$$d^{2n}_0(x,y)\le C \, d_0(x,y)^{\frac{n}{p}}d_f(x,y)^n\left(\mu_f(M)\right)^{-1}$$
We get therefore for all $p>1$,
$$d^{2n-\frac{n}{p}}_0(x,y)\le C\, d_f^n(x,y)\left(\mu_f(M)\right)^{-1},$$
which is the left-hand side inequality of Proposition \ref{prop:biHolder}.
\end{proof}

\section{Conformal deformations of $\mathbb R^n$ with scalar curvature in $L^{n/2}$}\label{sec:Rn}

In this section, we present some interesting application of our techniques to conformal deformations of $(\R^n, \eucl)$.

Let us first provide the following sufficient criterion to get strong $A_\infty$ weights.

\begin{thm}\label{theo:BHS-Rn} Let $f\colon \R^n\rightarrow \R$ be such that
$$\|d f\|_{L^{n}(\R^n, {\eucl})}\le \mathbf{I}$$ 
for some $\mathbf{I}>0$. Then the weight $e^{nf}$ is a strong $A_\infty$ weight with respect to ${\eucl}$, with constants $(\uptheta,\eta)$, where $\eta$ depends only on $g_0$ and $\uptheta$ depends only on $g_0$ and $\mathbf{I}$.
\end{thm}

This slightly improves Theorem 3.1 of \cite{BHS2}, since our statement is scale invariant.  The proof of this result is an immediate adaptation of our proof of Theorem \ref{theo:BHS}. We will not provide further details since it would be a mere repetition of the previous section.

Nevertheless, we show now that Theorem \ref{theo:BHS-Rn} implies Theorem \ref{theo:SAInfty-Rn}, which provides a completely new criterion to get Strong $A_\infty$ weights on $(\R^n, {\eucl})$.

\begin{proof}[Proof of Theorem \ref{theo:SAInfty-Rn}]

Let $g=e^{2f}\mathrm{eucl}$ be a conformal metric on the Euclidean space $\R^{n>2}$ such that 
\begin{enumerate}[i)]
\item $\vol(\R^n,g)=+\infty$;
\item $\int_{\R^n} |\scal_g|^{\frac n2}(x)\, \mathrm{dvol}_g(x)<+\infty$.
\end{enumerate}
We will show that $e^{nf}$ is a strong $A_\infty$ weight with respect to ${\eucl}$.

The following lemma shows that our infinite volume assumption is often satisfied.

\begin{lem}
A \emph{complete} conformal deformation $g=e^{2f}\mathrm{eucl}$ of the Euclidean metric on $\R^{n>2}$ that satisfies $\int_{\R^n} |\scal_g|^{\frac n2}(x)\, \mathrm{dvol}_g(x)<+\infty$ has necessary infinite volume.
\end{lem}
\begin{proof}
According to this hypothesis, the Riemannian manifold $(\R^n,g)$ satisfies the Sobolev inequality (\cite[Proposition 2.3]{CH} ): for some $C>0$,
$$\forall \varphi\in \cC^1_0(\R^n)\colon \left(\int |\varphi|^{\frac{2n}{n-2}}\ \mathrm{dvol}_g\right)^{1-\frac 2n} \le C \int |d\varphi|_g^{2}\ \mathrm{dvol}_g.$$
Hence $(\R^n,g)$ has an infinite volume, since the volume of geodesic balls is then bounded from below (\cite{Aku,CarronSMF}):
$$c(n,\mu) r^n \le \vol_g(B(x,r)).$$
\end{proof}

Now, our proof of Theorem \ref{theo:SAInfty-Rn} will be based on the following elementary lemma.

\begin{lem}\label{lem:Harmo-Constant}
Let $(M,g)$ be a complete Riemannian manifold that satisfies the elliptic Harnack inequality. Then harmonic functions which are bounded from below are constant.
\end{lem}

\begin{proof}
Harnack inequality establishes the existence of a constant $\upgamma>0$ such that for every geodesic ball $B(x,r)\subset M$ and every positive harmonic function $h$ defined over
$B(x,2r)$ one has
$$\sup_{y\in B(x,r)} h(y)\le \upgamma \inf_{y\in B(x,r)} h(y).$$
Hence if $h\colon M\rightarrow \R$ is bounded from below, setting $\uplambda=\inf_{y\in M} h(y)$, we get that for every geodesic ball
$$\sup_{y\in B(x,r)} (h(y)-\uplambda)\leq \gamma \inf_{y\in B(x,r)} (h(y)-\uplambda).$$
But as $r\to+\infty$, we have $$\lim_{r\to +\infty} \inf_{y\in B(x,r)} h(y)=\uplambda$$ therefore for all $y\in M$, $h(y) = \uplambda$.
\end{proof}

As an example that will be useful later, let us consider the case where $\bar g=e^{2u}\mathrm{eucl}$ is a conformal metric on the Euclidean space $\R^n$ and assume that 
$$\int_{\R^n} |du|^n(x)dx<+\infty.$$
It follows from Theorem \ref{theo:BHS-Rn} that $e^{nu}$ is a strong $A_\infty$ weight and that the metric $\bar g$ satisfies  the Poincar\'e inequality and the doubling condition: there are constants $\uptheta, \uplambda$ such that for all $\bar g$-geodesic balls $B(x,r)$, we have
$$\vol_{\bar g} B(x,2r)\le \uptheta \vol_{\bar g} B(x,r)$$
and, writing $\varphi_{B(x,r)}:=\fint_{B(x,r)} \varphi\, \mathrm{dvol}_{\bar g}$,
$$\forall \varphi\in \cC^1(B(x,2r))\colon \left\|\varphi-\varphi_{B(x,r)}\right\|_{L^2(B(x,r),{\bar g})}^2\le \left\|\nabla\varphi\right\|_{L^2(B(x,2r),{\bar g})}^2.$$

Moreover, by the results of Section \ref{ssec:SAInfty}, the metric space $(\R^n,\bar g,  \mathrm{dvol}_{\bar g})$ is Alhfors regu\-lar: there is a constant $\upeta$ such that for every $\bar g$-geodesic ball $B(x,r)$ :
$$\upeta r^n\le  \vol_{\bar g} B(x,r)\le r^n/\eta.$$
In particular by result in \cite{Gri, SC}, $\bar g$ satisfies the parabolic and the elliptic Harnack inequalities and if $n\ge 3$, then $(\R^n,\bar g)$ is a non parabolic manifold and has a positive Green kernel $G_{\bar g}(x,y)$ that satisfies the estimates 
\begin{equation}\label{eq:Green Control}
c d_{\bar g}^{2-n}(x,y)\le G_{\bar g}(x,y)\le C d_{\bar g}^{2-n}(x,y)
\end{equation}
or equivalently :
$$c \left(\vol_{\bar g}(\bB_{x,y})\right))^{-1+\frac 2n}\le G_{\bar g}(x,y)\le C  \left(\vol_{\bar g}(\bB_{x,y})\right))^{-1+\frac 2n}.$$

\medskip

Let us go back to our metric $g=e^{2f}\mathrm{eucl}$ with
$$
\vol(\R^n,g)=+\infty \quad \mbox{and} \quad \int_{\R^n} |\scal_g|^{\frac n2}(x)\, \mathrm{dvol}_g(x)<+\infty.
$$
Let $\Delta$ be the Euclidean Laplacian and let $A_n$ be the norm of the operator $d\Delta^{-1}\colon L^{\frac n2}\rightarrow L^n$. Since the $L^{n/2}$ norm of $\scal_g$ is finite, for $R$ large enough we have 
$$\frac{n-2}{4(n-1)}\left(\int_{\R^n\setminus \bB(R)}  |\scal_g|^{\frac n2}(x)\, \mathrm{dvol}_g(x)\right)^{\frac 2n}<\frac{1}{8A_n}.$$
Therefore, as was done in the proof of \pref{Gstate}, we can find $v\in L^{\frac n2}(\R^n,dx)$ such that
$$v-\left|d\Delta^{-1}v\right|^2=\frac{n-2}{4(n-1)} \scal_ge^{2f}\un_{\R^n\setminus \bB(R)}.$$

Now, let $u\in W^{2,\frac n2}_{\loc}(\R^n)$ be such that $\Delta u=v$, for instance
$u$ is a primitive of $d\Delta^{-1}v\in L^n$ or
$$u(x)=c_n\int_{\R^n}\left(\frac{1}{\|x-y\|^{n-2}}-\frac{1}{\|y\|^{n-2}}\right)v(y)dy.$$
We have 
$$\Delta e^u=\frac{n-2}{4(n-1)} \scal_ge^{2f}\un_{\R^n\setminus \bB(R)}e^u.$$ 
Writing $\upphi=e^{\frac{n-2}{2}f}$, Yamabe equation (\ref{Yamabe}) implies that
$$\Delta\upphi=\frac{n-2}{4(n-1)} \scal_ge^{2f}\upphi.$$
Hence if we define $\uppsi=\upphi e^{-u}$ we get
$$\Delta\uppsi-\la du,d\uppsi\ra=\frac{n-2}{4(n-1)} \scal_ge^{2f}\un_{ \bB(R)}\uppsi.$$
As $du\in L^n$, the metric $\bar g=e^{\frac{4}{n-2}u}\mathrm{eucl}$  induces a strong $A_\infty$ weight, but
$$\Delta_{\bar g} \uppsi=e^{-\frac{4}{n-2}u}\left(\Delta\uppsi-\la du,d\uppsi\ra\right)=\frac{n-2}{4(n-1)} \scal_ge^{2f}\un_{ \bB(R)}e^{-\frac{4}{n-2}u}\uppsi.$$
As $\Delta_{\bar g} \uppsi$ has compact support and is in $L^p$ for every $1\le p<+\infty$, we can find a solution of the equation
$$\Delta_{\bar g} \upxi= \Delta_{\bar g} \uppsi,$$ for instance let us choose
$$\upxi(x)=\int_{\R^n} G_{\bar g}(x,y)\Delta_{\bar g} \uppsi(y)\mathrm{dvol}_{\bar g}(y).$$

Since $\supp\Delta_{\bar g} \uppsi\subset K$, it follows from (\ref{eq:Green Control}) that there exists $C>0$ such that for all $x\in M\backslash K$,
$$\upxi(x)\le C\left(\int_{K} G_{\bar g}^q(x,y)\mathrm{dvol}_{\bar g}(y)\right)^{\frac 1q}\le C d_{\bar g}(x,K)^{2-n}.$$
In particular, $\int_{\R^n} \left| \upxi(y)\right|^{\frac{2n}{n-2}}\mathrm{dvol}_{\bar g}(y)<+\infty$,
the map $\upxi$ tends to zero at infinity and is in $W^{2,p}_{\loc}(\R^n)$ for every $1\le p<+\infty$. Therefore $\upxi$ is a bounded function and 
$\uppsi-\upxi$ is a $\bar g$-harmonic function that is bounded from below: it is hence a constant function.
If this constant is $0$ then $\uppsi=\upxi$ and 
$$\vol(\R^n,g) = \int_{\R^n} \uppsi(y)^{\frac{2n}{n-2}}\mathrm{dvol}_{\bar g}(y)=\int_{\R^n} \upxi^{\frac{2n}{n-2}}(y)dy<+\infty.$$
This is contradictory with our hypothesis hence there is a positive constant $c$ such that
$$\uppsi=c+\upphi,$$ hence there are positive constants $\alpha,A$ such that
$$\alpha\le \uppsi\le A.$$
This implies that the metric $g=\upphi^{\frac{4}{n-2}}\mathrm{eucl}$  and $\bar g=e^{2u}\mathrm{eucl}$ are bi-Lipischitz, hence since $e^{nu}$ is a strong $A_\infty$ weight with respect to $g_{\eucl}$, $e^{nf} = \upphi^{\frac{2n}{n-2}}$ is also.

\end{proof}

We emphasize again on the fact that this implies very strong geometric constraints on conformal deformations of $(\R^n, g_{\eucl})$ with $\displaystyle \int_{\R^n} |\scal_g|^{\frac n2}(x)\, \mathrm{dvol}_g(x)<+\infty$ and $\vol_g(\R^n,g)=+\infty$. It implies in particular Poincar\'e inequality, Euclidean-type isoperimetric inequality and Sobolev inequalities, which were unknown on such examples.

\section{Gromov-Hausdorff limit and conformal distances}\label{sec:LimDist}

Let us conclude by presenting a natural question about the limiting distances which arise in our Theorem \ref{theo:SAInfty}, which we could not answer. 

\medskip

Let us first recall that to any strong $A_\infty$ weight $e^{n f_\infty}$ (with respect to $d{\rm vol}_{g_0}$) is naturally associated a distance $d_f$, defined as follows. For all $x,y\in M$, we write agin $\bB_{xy}$ the largest $g_0$-geodesic ball whose diameter is the $g_0$-geodesic segment $[xy]$. Let us write 
$$\delta_f(x,y) = \left(\frac 1 {\omega_n}\mu_{f}(\bB_{xy})\right)^{1/n} = \left(\frac 1 {\omega_n}\int_{\bB_{xy}}e^{n f(x)}d{\rm vol}_{g_0}(x)\right)^{1/n},$$ 
where $\omega_n$ is the Euclidean volume of the unit $n$-ball, and define
\begin{equation}\label{eq:DistSemmes}
d_{f}(x,y) = \lim_{\epsilon\to 0} \inf \left\{ \sum_{i = 0}^N \delta_f(x_i, x_{i+1}) \; ; \; x_0 = x, x_{N+1} = y , d_0(x_i, x_{i+1})\leq \epsilon \right\}.
\end{equation}
It can be easily shown that $e^{n f}\in L^1(M)$ is a strong $A_\infty$ weight if and only if $d_f$ is a distance (see \cite{Semmes1} for more details). Moreover, an elementary approximation argument shows that if $f$ is continuous, then $d_f$ coincides with the classical Riemannian distance associated to the Riemannian metric $g_f = e^{2f} g_0$.

\medskip

Let us now consider a sequence $(g_k = e^{2 f_k} g_0)_{k\in \N}$ of conformal metrics satisfying the hypotheses of Corollary \ref{coro:biHolder}. Up to extracting a subsequence, the sequence of associated metric space $(M, d_{g_k})$ converges in Gromov-Haudsorff topology to some distance $d_\infty$ on $M$, which is bi-H\"older equivalent to $d_0$. Moreover, still up to extracting a subsequence, the Riemannian volume measure $\mu_k = d{\rm vol}_{g_k} = e^{n f_k}d{\rm vol}_{g_0}$ converge in the weak$-*$ topology to a positive measure $\mu_\infty$ on $M$. Since all the $\mu_k$ are uniformly strong $A_\infty$ weights with respects to $d{\rm vol}_{g_0}$, $\mu_\infty$ is also: there exists hence $e^{n f_\infty}\in L^1(M)$ which is still a strong $A_\infty$ weight with respect to $g_0$. This weight is therefore associated to a distance $d_{f_\infty}$ defined by (\ref{eq:DistSemmes}). \emph{Does $d_\infty$ and $d_{f_\infty}$ coincide ?} 

\medskip
The reader will can check that in the Example 2 presented in our introduction, where the conformal factors converge in $W^{2, \frac n 2}$ to the singular conformal factor given by \ref{eq:SchrodEx}) and the distances converge in the Gromov-Hausdorff topology, these distances indeed coincide.

\medskip
One can easily show from semi-continuity arguments that $d_\infty \leq d_{f_\infty}$. Let us present now a striking example of pathological behaviour which shows that in general, for a sequence of (unifom) strong $A_\infty$ weights $(e^{n f_k})_{k\in \N}$, the converse inequality is not true. The following examples of conformal deformations is due to D. Burago, in \cite{Burago}

Let us consider the standard flat torus $\ds \left(\bT^n=\R^n/\Gamma, \mathrm{eucl}\right)$, where $\Gamma=\left(2\pi \Z\right)^n$. For all $\ell>1$, we consider the conformal deformation of the euclidean metric $g_0$ given by $g_\ell = e^{2 f_\ell} \mathrm{eucl}$, where
$$
e^{nf_\ell}(x_1,\dots,x_n) =1-\frac12 \cos(\ell x_1).
$$
Since $e^{nf_\ell}$ is uniformly bounded from above and below, all the metric space $(\bT^n,d_\ell=d_{f_\ell})$ are uniformly bilipschitz to $\bT^n$ equiped with the Euclidean distance and the $(e^{nf_\ell})_{\ell\geq 1}$ are uniformly strong $A_\infty$ weights with respect to the Lebesgue metric. Nevertheless, the distance associated to the limit of the volume measures and the limit of distances are radically distinct as shows the following result.

\begin{thm}
Under the above notations, as $\ell\to +\infty$,
\begin{itemize}
\item the sequence $(e^{n f_\ell})_\ell$ converges weakly in $L^1(\bT)$ to the function $1$;
\item the sequence of metric spaces $(\bT^n,d_\ell), \ell \in \N$ converges in the Gromov-Hausdorff topology to metric space
$(\bT^n,d_\infty)$ where $d_\infty$ is the Finsler metric associated to the stable norm of the periodic metric on $\R^n$ given by
$$\bar g=\left(1-\frac12 \cos(x_1)\right)^{\frac2n} \mathrm{eucl}.$$
\end{itemize}
\end{thm}

\begin{proof}
The first assertion is classical. As for the second assertion, let us recall that if $\bar d$ is the Riemannian distance on $\R^n$ associated to the metric $\bar g$ then the associated stable norm $\|\bullet\|_*$ is
defined for all $x\in \R^n$ by 
$\|x\|_*=\lim_{t\to+\infty}\frac{1}{t} \bar d(0,tx).$
Using the change of variable $\xi_i=\ell x_i, i=1\dots n$, we get that $(\bT^n,d_\ell=d_{f_\ell})$ is isometric to the quotient
$\R^n/(\ell \Gamma)$ with the metric $\frac{1}{\ell^2}\bar g$.
More exactly let $\pi\colon\R^n\rightarrow \bT^n=\R^n/\Gamma$ the quotient map. if $x,y\in \R^n$ then
$$d_\ell\left(\pi(x),\pi(y)\right)=\inf_{\gamma\in \Gamma}\frac{1}{\ell} d(x,y+\ell\gamma)$$ 
If $d_\infty\left(\pi(x),\pi(y)\right)=\inf_{\gamma\in \ell\Gamma}\|x-y-\gamma\|_*$, then we clearly have
$$\sup_{x,y\in \R^n} \left| d_\ell\left(\pi(x),\pi(y)\right)-d_\infty\left(\pi(x),\pi(y)\right)\right|\le \frac{C}{\ell}.$$
Hence the result.
\end{proof}

We emphasize on the fact that in this example, $d_{f_\infty} = d_{eucl}$, whereas $d_\infty$ is a Finsler non-Riemannian metric: it cannot in any weak sense be considered as a conformal metric for some limiting conformal factor.

\medskip

The reader will easily notice than in the previous example, the volume of the conformal metrics are uniformly bounded, but the $L^{n/2}$-norm of the scalar curvature blows up. Under the $L^{n/2}$ assumptions on the scalar curvature which ensures Theorem \ref{theo:SAInfty} and Corollary \ref{coro:biHolder}, we have a stronger control on the geometrical and analytical behaviour of the conformal factors than the information given by the strong $A_\infty$ bounds. Nevertheless, we cannot show yet equality for $d_\infty$ and $d_{f_\infty}$ under these hypotheses, nor give a counterexample.

\end{document}